\documentclass{amsart}

\usepackage{color}
\usepackage{graphicx}%
\usepackage{amsmath,amssymb,amsfonts,amsthm,amsxtra}%
\usepackage{mathrsfs}%
\usepackage{xcolor}%
\usepackage{bm}
\usepackage{mathtools}
\mathtoolsset{showonlyrefs}
\usepackage[bookmarks=false,unicode,linktocpage=true,pagebackref=false,hyperindex=true,pageanchor,colorlinks=true,linkcolor=blue,linkbordercolor=red,pdfborderstyle={/S/U/W 1}]{hyperref}
\usepackage{dsfont}
\usepackage{bm}

\newtheorem{theorem}{Theorem}
\newtheorem{proposition}{Proposition}% 
\newtheorem{lemma}{Lemma}% 
\newtheorem{remark}{Remark}%
\newtheorem{definition}{Definition}%

\definecolor{magenta}{rgb}{1,0,1}
\definecolor{cyan}{rgb}{0,1,1}
\definecolor{obsgrey}{rgb}{0.7,0.7,0.7}
\definecolor{orange}{rgb}{1,0.5,0}
\definecolor{brown}{rgb}{0.5,0.25,0}
\definecolor{violet}{rgb}{0.25,0.0,0.65}
\definecolor{darkblue}{rgb}{0,0,0.5}
\definecolor{darkgreen}{rgb}{0,0.5,0}
\definecolor{darkred}{rgb}{0.65,0,0}
\definecolor{lightgrey}{rgb}{0.6,0.6,0.6}
\definecolor{vanish}{rgb}{0.6,0.6,0.6}
\definecolor{lightcyan}{rgb}{0.3,1,1}
\definecolor{medcyan}{rgb}{0,0.7,0.7}
\definecolor{result}{rgb}{0.5,0,0}
\definecolor{lightviolet}{rgb}{0.6,0.0,0.9}
\definecolor{modified}{rgb}{0.3,0.3,1}
\definecolor{midgreen}{rgb}{0,0.8,0}

\newcommand{\mysection}[1]{\section{#1}}
\newcommand{\myref}[1]{\ref{#1}}
\newcommand{\mylabel}[1]{\label{#1}}
\newcommand{\myfiglabel}[1]{\label{#1}}
\newcommand{\myeqlabel}[1]{\label{#1}}
\newcommand{\myeqref}[1]{\eqref{#1}}

\newcommand{\union}{\cup}
\newcommand{\closure}[1]{\overline{#1}}
\newcommand{\setdiff}{\backslash}

\renewcommand{\Re}{\operatorname{Re}}
\renewcommand{\Im}{\operatorname{Im}}

\newcommand{\cartp}{\times}
\newcommand{\boi}[2]{{]#1,#2[}}

\newcommand{\roi}[2]{{[#1,#2[}}
\newcommand{\cli}[2]{{[#1,#2]}}

\newcommand{\sign}{\operatorname{sign}}

\newcommand{\Det}{\operatorname{det}}
\newcommand{\adjt}{\operatorname{adj}}

\newcommand{\eqv}{\Leftrightarrow}
\newcommand{\impl}{\Rightarrow}
\newcommand{\qeq}{\quad\eqv\quad}

\newcommand{\NN}{\mathds{N}}
\newcommand{\RR}{\mathds{R}}
\newcommand{\ZZ}{\mathds{Z}}

\newcommand{\bdry}{\partial}

\newcommand{\pd}[1]{\partial_{#1}}

\newcommand{\conv}{\rightarrow}

\newcommand{\upconv}{\nearrow}
\newcommand{\dnconv}{\searrow}
\newcommand{\defm}[1]{\emph{#1}}
\newcommand{\supeq}[2]{\mathord{\overbrace{\mathop{#1}}^{#2}}}
\newcommand{\subeq}[2]{\mathord{\underbrace{\mathop{#1}}_{#2}}}
\newcommand{\csep}{\quad,\quad}
\newcommand{\half}{{\frac12}}

\newcommand{\dotp}{\cdot}

\newcommand{\trace}{\operatorname{trace}}

\newcommand{\Lap}{\Delta}
\newcommand{\Lapd}{\Delta_{\dpq}}

\newcommand{\eqs}{\sim}
\newcommand{\topref}[2]{\overset{\text{\myeqref{#1}}}{#2}}
\newcommand{\toprefb}[3]{\overset{\text{\myeqref{#1}}}{\underset{\text{\myeqref{#2}}}{#3}}}
\newcommand{\sett}[1]{\{#1\}} 
\newcommand{\selset}[2]{\{#1~:~#2\}}
\newcommand{\spC}{\mathcal{C}}
\newcommand{\Ck}[1]{\spC^{#1}}

\newcommand{\Conce}{\Ck1}

\newcommand{\isect}{\cap}
\newcommand{\Cinf}{\Ck\infty}

\newcommand{\find}{\operatorname{ind}}

\newcommand{\const}{\text{const}}

\newcommand{\abs}[1]{\vert#1\vert}   
\newcommand{\vlen}[1]{\vert#1\vert}
\newcommand{\vlenq}[1]{\vlen{#1}^2}

\newcommand{\eps}{\epsilon}

\newcommand{\rcori}{\rad_\cori}
\newcommand{\rw}{\rad_w}
\newcommand{\ro}{\rad_o}
\newcommand{\ri}{\rad_i}

\newcommand{\pqvar}{p} 
\newcommand{\dpqvar}{\tilde\pqvar} 
\newcommand{\pq}{{\vec\pqvar}} 
\newcommand{\pqin}{\pqi} 
\newcommand{\pqi}{\pq_i} 
\newcommand{\pqw}{\pq_w} 
\newcommand{\pqwx}{\pqvar_w^1} 
 
\newcommand{\pqix}{\pqvar_i^1} 
\newcommand{\pqiy}{\pqvar_i^2} 
\newcommand{\pqs}{\pq_s}

\newcommand{\pqo}{\veczero} 
\newcommand{\pqcori}{\pq_\cori}
\newcommand{\dpqcori}{\dpq_\cori}
\newcommand{\Dompq}{\Dom_{\pq}}
\newcommand{\cDompq}{\closure\Dom_{\pq}}
\newcommand{\dDompq}{\Dom_{\dpq}}
\newcommand{\cdDompq}{\closure\Dom_{\dpq}}
\newcommand{\Axispq}{\Axis_{\pq}}
\newcommand{\cAxispq}{\closure\Axis_{\pq}}
\newcommand{\dAxispq}{\Axis_{\dpq}}

\newcommand{\Wedgepq}{\Wedge_{\pq}}
\newcommand{\cWedgepq}{\closure\Wedge_{\pq}}
\newcommand{\dWedgepq}{\Wedge_{\dpq}}

\newcommand{\Shockpq}{\Shock_{\pq}}
\newcommand{\cShockpq}{\closure\Shock_{\pq}}
\newcommand{\dShockpq}{\Shock_{\dpq}}

\newcommand{\pqpot}{\Psi}
\newcommand{\pqpotshift}{\pqpot_*}

\newcommand{\nabp}{\nabla_\pq}
\newcommand{\ttpsym}{t}
\newcommand{\nnpsym}{n}
\newcommand{\ttp}{\vec\ttpsym}
\newcommand{\nnp}{\vec\nnpsym}
\newcommand{\curv}{\kappa}

\newcommand{\dpq}{{\vec{\dpqvar}}}
\newcommand{\dpqw}{\dpq_w}

\newcommand{\dpqi}{\dpq_i}
\newcommand{\dpqo}{\dpq_o}

\newcommand{\defeq}{:=}
\newcommand{\eqdef}{=:}

\renewcommand{\vec}[1]{\bm{#1}}
\newcommand{\veczero}{\vec{0}}
\newcommand{\omegav}{{\omega}} 
\newcommand{\cori}{{\tau}}

\newcommand{\omegalo}{\omegav_0}
\newcommand{\omegahi}{\omegav_1}
\newcommand{\omegaa}[1]{\omegav_{#1}}
\newcommand{\omegazero}{\omegaa0}
\newcommand{\omegaone}{\omegaa1}
\newcommand{\jside}{{j}}
\newcommand{\omegajside}{\omegaa{\jside}}

\newcommand{\vsym}{u}

\newcommand{\vv}{\vec\vsym}
\newcommand{\vx}{\vsym^\xe}
\newcommand{\vy}{\vsym^\ye}
\newcommand{\vvn}{\vv_n}

\newcommand{\vvs}{\vv_s}

\newcommand{\vvw}{\vv_w}
\newcommand{\vvwmir}{\bar\vv^w}
\newcommand{\vvsmir}{\bar\vv^s}

\newcommand{\vvc}{\vv_c}
\newcommand{\vvca}{\vv_{c*}}

\newcommand{\jmp}[1]{[#1]}

\newcommand{\jv}{j}
\newcommand{\jj}{\vec{\jv}}
\newcommand{\jmax}{\jv_\max}
\newcommand{\jjs}{\jj_s}
\newcommand{\jjw}{\jj_w}
\newcommand{\jjc}{\jj_c}
\newcommand{\jji}{\jj_i}
\newcommand{\jjn}{\jj_n}
\newcommand{\jx}{\jv^x}
\newcommand{\jy}{\jv^y}

\newcommand{\jyw}{\jv^y_w}

\newcommand{\firstcoeff}{\vec{a}}
\newcommand{\lpr}{\ell}
\newcommand{\allop}{\mathcal{T}} 
\newcommand{\inop}{\mathcal{L}}
\newcommand{\bdop}{\mathcal{B}}

\newcommand{\bt}{\beta} 
\newcommand{\btvec}{\vec\bt}

\newcommand{\domsob}[2]{\mathcal{X}^{#1}_{#2}}
\newcommand{\extdomsob}[2]{\hat{\mathcal{X}}^{#1}_{#2}}
\newcommand{\ransob}[2]{\mathcal{Y}^{#1}_{#2}}
\newcommand{\dvlb}[2]{\mathscr{D}^{#1}_{#2}V}
\newcommand{\rvlb}[2]{\mathscr{R}^{#1}_{#2}V}

\newcommand{\sobl}{l}

\newcommand{\mat}[1]{\mathbf{#1}}
\newcommand{\Acosym}{A} 
\newcommand{\Aco}{\mat{\Acosym}}

\newcommand{\arccot}{\operatorname{arccot}}

\newcommand{\nuia}{\mu}
\newcommand{\nuiasonic}{\nuia_i^1} 
\newcommand{\nuiacrit}{\nuia_i^c} 
\newcommand{\nuiavan}{\nuia_i^v} 
\newcommand{\turnap}{\theta}
\newcommand{\turnai}{\turnap_w}
\newcommand{\turnamax}{\turnap_c}

\newcommand{\dexsym}{e}
\newcommand{\dex}{\vec{\dexsym}}

\newcommand{\bbsym}{b}
\newcommand{\bzero}{\bbsym^0}
\newcommand{\dbb}{\vec\bbsym}

\renewcommand{\min}{{\operatorname{min}}}
\renewcommand{\max}{{\operatorname{max}}}

\newcommand{\densstag}{\dens_0}
\newcommand{\csndstag}{\csnd_0}

\newcommand{\dbcang}{\beta}
\newcommand{\dbcanga}[1]{\dbcang_{#1}}
\newcommand{\dbcangjside}{\dbcanga\jside}
\newcommand{\dbcangow}{\dbcanga1^w}
\newcommand{\omegaow}{\omegaone}

\newcommand{\wtvec}{\vec{\wt}}
\newcommand{\wt}{\zeta}
\newcommand{\wtw}{\wt_w}
\newcommand{\wtoi}{\wt_{oi}}

\newcommand{\wtcori}{\wt_\cori}
\newcommand{\critex}{{\alpha}}
\newcommand{\ice}{k}
\newcommand{\sigk}{\varsigma}

\newcommand{\pr}{\pd\rad}
\newcommand{\pom}{\pd\omegav}
\newcommand{\radpr}{\rad\pr}

\newcommand{\hh}{h}
\newcommand{\hw}{\hh^\Wedge}
\newcommand{\hs}{\hh^\Shock}

\newcommand{\xx}{{\vec x}}
\newcommand{\xe}{{x}}
\newcommand{\ye}{{y}}
\newcommand{\yew}{\ye^w}

\newcommand{\cutfn}{\chi}
\newcommand{\cutcori}{\cutfn_\cori}
\newcommand{\cutfni}{\cutfn_i}
\newcommand{\cutfno}{\cutfn_o}
\newcommand{\cutfnw}{\cutfn_w}
\newcommand{\rad}{r}

\newcommand{\Wedge}{W}

\newcommand{\Axis}{X}
\newcommand{\Shock}{S}
\newcommand{\Dom}{\Omega}

\newcommand{\nzplDom}{G}

\newcommand{\Wedgexy}{\Wedge}
\newcommand{\Shockxy}{\Shock}
\newcommand{\Domxy}{\Dom}
\newcommand{\Axisxy}{\Axis}
\newcommand{\cDomxy}{\closure\Dom}

\newcommand{\cShockxy}{\closure\Shock}
\newcommand{\cAxisxy}{\closure\Axis}

\newcommand{\eWedgexy}{\hat\Wedge}
\newcommand{\eShockxy}{\hat\Shock}
\newcommand{\eDomxy}{\hat\Dom}

\newcommand{\ecDomxy}{\closure{\hat\Dom}}

\newcommand{\dts}{\vec{t}}

\newcommand{\shocknsym}{n}

\newcommand{\shockn}{{\vec\shocknsym}}

\newcommand{\wallnsym}{n}

\newcommand{\walln}{\vec\wallnsym}

\newcommand{\idm}{\mat I} 

\newcommand{\xxcS}{\xx_c^\Shock}
\newcommand{\xxcW}{\xx_c^\Wedge}
\newcommand{\xxn}{\xx_n}
\newcommand{\xxo}{\veczero}
\newcommand{\xxw}{\xx_w}

\newcommand{\Mach}{M}

\newcommand{\densi}{\dens_i}
\newcommand{\csndi}{\csnd_i}

\newcommand{\vvi}{\vv_i}
\newcommand{\vxi}{\vsym_i^\xe}

\newcommand{\vpot}{\phi}
\newcommand{\genphi}{\varphi}

\newcommand{\vpoti}{\vpot_{i}}
\newcommand{\stf}{\psi}
\newcommand{\stfi}{\stf_i}
\newcommand{\substf}{\underline\psi}

\newcommand{\vang}{\turnap}
\renewcommand{\Mach}{M}

\newcommand{\Machi}{\Mach_i}

\newcommand{\pressure}{P}
\newcommand{\pressurei}{\pressure_i}

\newcommand{\dens}{\varrho}

\newcommand{\gisen}{\gamma}
\newcommand{\nabx}{\nabla_{\xx}}

\newcommand{\nabd}{\nabla_{\dpq}}

\newcommand{\spm}{s}
\newcommand{\spmi}{s_i}
\newcommand{\hpm}{h}
\newcommand{\hpmi}{\hpm_i}
\newcommand{\Hpm}{H}
\newcommand{\Hpmi}{\Hpm_i} 

\newcommand{\csnd}{c}

\newcommand{\mytitle}{Weak-type bow shocks}

\newcommand{\myabstract}{%
  An existence proof is given for bow shocks ahead of blunt convex wedges that are transonic and of weak type at infinity, including cases for arbitrarily small shock strength.
  It is also proven that any such flow is necessarily univalent, meaning the position-to-velocity map is injective. 
}
\title{\mytitle}
\author{Volker W. Elling}
\begin{document}
\maketitle
\begin{abstract}\myabstract\end{abstract}

\mysection{Introduction}

Consider supersonic flow of ideal gas onto a blunt smooth convex solid wedge that is symmetric across the flow axis and has straight sides in a neighbourhood of infinity.
If the straight sides are extended to their meeting point (tip) they form a sharp wedge (fig.\ \myref{fig:sharp-wedge}).
If the wedge half-angle $\turnai$ is sufficiently small, the expected flow at a sharp wedge has attached shocks at the wedge tip,
with each one of two possible choices called weak-type and strong-type. 

If $\turnai$ increases to the \defm{critical angle} $\turnamax$ the two types converge to a single critical-type shock; for larger $\turnai$ there is no solution.
If $\turnai$ decreases to zero, the strong-type shock approaches the normal shock perpendicular to $\vvi$ with post-shock velocity $\vv$ parallel to $\vvi$ but strictly smaller;
the strength of the weak-type shock converges to zero while its slope converges not to the wedge surface but to an acoustic wave with \defm{Mach angle} $\arcsin(1/\Machi)$ to $\vvi$,
where $\Machi$ is the Mach number in the constant pre-shock region.

Sharp wedges permit attached shocks; the related problems have been studied extensively,
both for purely subsonic and for purely supersonic post-shock region (e.g.\ \cite{yongqian-zhang-wedge-glimm,bae-liao-2026,blokhin-tkachev-baldan,bae-review-sonic,qin-wang-kyungwoo-song-convex-cornered,blokhin-tkachev-yegitov,li-xu-yin-threedweakstrong,li-witt-yin-2014,hu-qu-2026} and references therein).
For a \emph{blunt} wedge the shock is observed to be a detached \defm{bow shock} ahead of the body (fig.\ \myref{fig:bluntweakstrong}).
Although bow shocks are ubiquitous in applications ranging from engineering to astrophysics, mathematical existence proofs for these flows are hard to find.

\begin{figure}
  \centerline{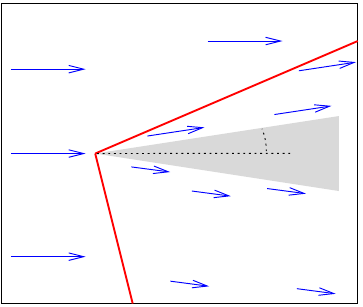}
  \caption{Supersonic flow onto a symmetric sharp wedge}
  \myfiglabel{fig:sharp-wedge}
\end{figure}

\begin{figure}
  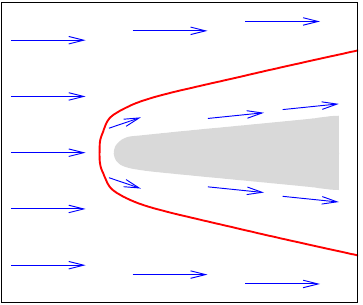\hfil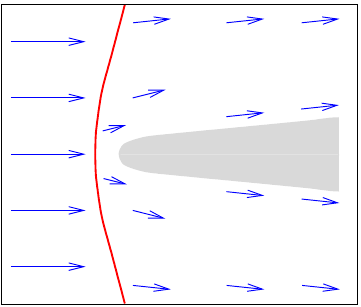%
  \caption{Left: blunt convex symmetric wedge with bow shock converging at infinity to weak-type. Right: another wedge with convergence to strong-type.}
  \myfiglabel{fig:bluntweakstrong}
\end{figure}

Behind detached shocks flow is always subsonic at some point of the body,
but often supersonic near infinity,
separated by shocks or sonic transitions that pose some interesting technical challenges.
But even for purely subsonic post-shock region there are few rigorous proofs. 
Some prior work can be found in \cite{chen-christoforou-jegdic} for bow shocks ahead of a sharp wedge,
and \cite{bae-xiang-blunt} on the $\Machi\conv\infty$ limit at blunt wedges with strong-type shock at infinity.

In this article we develop new techniques for proving the first existence results for \emph{weak} type for some convex wedges and some $\Machi$, with uniformly subsonic post-shock region.
This is important because in nature and numerics the commonly observed shock is the weak-type one.
As von Neumann \cite[p.\ 254]{neumann-1943} states:
``Now the [strong-type] possibility amounts to an infinitesimal shock producing a reflected shock of finite strength.
This runs counter to common sense, even for energetic reasons. Besides we know from experiments that sound waves are reflected as sound waves.''
Although this remark refers to regular reflection of an incident shock of infinitesimal strength,
it can equally be applied to a wedge with infinitesimal angle $\turnai$.
Courant-Friedrichs \cite[p.\ 317]{courant-friedrichs} consider ramps rather than wedges,
noting that ``it has been frequently stated that the [strong-type] is unstable and that, therefore, only the [weak-type] could occur'',
although they propose that the strong type could appear ``under appropriate circumstances''.
\cite[section 6]{elling-massmodel} provides a simple heuristic mass-flux model clarifying the nature of this instability for shocks attached to sharp wedges
(see also \cite{elling-sonic-potf,elling-detachment,teshukov,elling-liu-rims05}).

After recalling some well-known background we state the main results in section \myref{section:mainresults}, with a brief overview of the technical obstacles and rest of the paper.

\mysection{Equations}\mylabel{section:equations}

The equations and their derivation are well-known, so we go over them as quickly as establishing the notation and assumptions permits. 

Consider steady planar inviscid irrotational flow. Pressure is determined by a $\gisen$-law equation of state
\begin{alignat}{9} \pressure = \pressurei(\dens/\densi)^\gisen \myeqlabel{eq:polytropic-pressure}\end{alignat}
where $\pressurei,\densi$ are positive constants fixed by a constant pre-shock flow, with \defm{isentropic coefficient} $\gisen\in\boi1\infty$. 
(Our result could be extended to $\pressure$ that are more general analytic functions of $\dens$, under conditions explained in \cite{elling-polarpotf},
but the polytropic case is by far the most important.)
We have assumed that total enthalpy $\Hpm=\hpm+\half\vlenq{\vv}$ and entropy per mass $\spm$ are global constants $\Hpmi,\spmi$, which is necessary for irrotationality in the form
\begin{alignat}{9} 0 = \nabx^\perp\dotp\vv \myeqlabel{eq:irrotationality}\end{alignat}
(where $(a,b)^\perp=(-b,a)$ is rotation by 90 degree counterclockwise); 
these assumptions combined imply (at least locally for smooth flow) conservation of momentum and energy.
$\hpm$ is enthalpy per mass, here
\begin{alignat}{9} \hpm = \hpmi (\dens/\densi)^{\gisen-1} \csep \hpmi = \frac{\gisen\pressurei}{(\gamma-1)\densi}, \myeqlabel{eq:hpm-dens}\end{alignat}
satisfying the usual relation
\begin{alignat}{9} \frac{\partial\hpm}{\partial\dens}_{|\spm} = \frac{\csnd^2}{\dens} \myeqlabel{eq:dhpm-ddens} \end{alignat}
where subscripts $|...$ indicate which coordinates are held constant for taking partial derivatives; 
\begin{alignat}{9} \csnd^2 = \frac{\partial\pressure}{\partial\dens}_{|\spm} > 0 \myeqlabel{eq:csnd}\end{alignat}
with $\csnd>0$ \defm{speed of sound}.

\myeqref{eq:irrotationality} implies (on simply connected domains)
\begin{alignat}{9} \vv=\nabx\vpot  \myeqlabel{eq:vpot} \end{alignat}
for a scalar \defm{velocity potential} $\vpot$. The Bernoulli relation in the form
\begin{alignat}{9}  \text{constant $\Hpmi$} = \Hpm = \half\vlenq{\vv} + \hpm  \myeqlabel{eq:bern-vv}  \end{alignat}
combined with \myeqref{eq:hpm-dens} determines $\dens$ as an analytic function of $\vv=\nabx\vpot$ by the analytic implicit function theorem (\cite[thm.\ 2.3.1]{krantz-parks}), using
\begin{alignat}{9}  \frac{\partial\Hpm}{\partial\dens}_{|\spm,\vv}  &=  \frac{\partial\hpm}{\partial\dens}_{|\spm}\topref{eq:dhpm-ddens}{=} \csnd^2\dens^{-1} > 0  ,
\\  \frac{\partial\Hpm}{\partial\vv}_{|\spm,\dens}  &=  \vv^T , \end{alignat}
showing
\begin{alignat}{9} \frac{\partial\dens}{\partial\vv}_{|\spm,\Hpm}  &=
- \Big( \frac{\partial\Hpm}{\partial\dens}_{|\spm,\vv} \Big)^{-1}  \frac{\partial\Hpm}{\partial\vv}_{|\spm,\dens}  =
-\csnd^{-2}\dens\vv^T . \myeqlabel{eq:ddens-dvv} \end{alignat}
So the continuity equation
\begin{alignat}{9} 0 = \nabx\dotp \jj \csep \jj=\dens\vv \myeqlabel{eq:continuity}\end{alignat}
for \defm{mass flux} $\jj$ becomes a second-order scalar equation for $\vpot$. Assuming smoothness we can pass to the non-divergence form
\begin{alignat}{9} 0 = \Aco:\nabx^2\vpot , \myeqlabel{eq:vpot-interiorpde-xy} \end{alignat}
with $\mat{B}:\mat{C}\defeq\trace(\mat{B}^T\mat{C})$ and $\nabla^2=\nabla\nabla^T$ (like $\vec a^2=\vec a\vec a^T$), where we calculated
\begin{alignat}{9} \frac{\partial(\dens\vv)}{\partial\vv}_{|\spm,\Hpm} &= \dens\idm + \vv\frac{\partial\dens}{\partial\vv}_{|\spm,\Hpm}
\\&\topref{eq:ddens-dvv}{=} \dens\big(\subeq{\idm-(\frac{\vv}{\csnd})^2}{\eqdef\Aco}\big). \myeqlabel{eq:ddensvv-dvv}\end{alignat}

$\Aco$ has eigenvectors $\vv$ for eigenvalue $1$ and $\vv^\perp$ for eigenvalue $1-\Mach^2$, for \defm{Mach number}
\begin{alignat}{9} \Mach = \frac{\vlen\vv}{\csnd}, \end{alignat}
so the 2nd order PDE \myeqref{eq:vpot-interiorpde-xy} for scalar $\vpot$ is elliptic exactly when $\Mach<1$ (subsonic flow), hyperbolic for $\Mach>1$ (supersonic). 

We will need an equivalent formulation: \myeqref{eq:continuity} implies
\begin{alignat}{9} \jj = -\nabx^\perp\stf \end{alignat}
for a scalar \defm{stream function} $\stf$ (at least on simply connected domains, or in presence of helpful boundary conditions). The Bernoulli relation in the form
\begin{alignat}{9} \text{constant $\Hpmi$} = \Hpm = \hpm + \half \dens^{-2} \vlenq{\jj} \myeqlabel{eq:bernj}\end{alignat} 
defines by \myeqref{eq:hpm-dens} $\dens$ as an analytic function of $\jj$; for applying the implicit function theorem we calculate
\begin{alignat}{9} \frac{\partial\Hpm}{\partial\jj}_{|\spm,\dens} &= \dens^{-2}\jj^T , 
\\ \frac{\partial\Hpm}{\partial\dens}_{|\spm,\jj} &\topref{eq:dhpm-ddens}{=} \frac{\csnd^2}{\dens} - \frac{\vlenq{\jj}}{\dens^3} = \frac{\csnd^2}{\dens} ( 1-\Mach^2 ).
\myeqlabel{eq:dHpm-ddens-atjj}\end{alignat}
For subsonic flows \myeqref{eq:dHpm-ddens-atjj} is positive, producing a map $\jj\mapsto\dens$ defined on the closure of the disk
\begin{alignat}{9} \selset{ \jj\in\RR^2}{ \vlen{\jj} < \jmax }, \myeqlabel{eq:jdisk}\end{alignat}
where $\jmax$ is the \defm{maximal mass flux}, a function of the constants $\Hpmi,\gisen,\densi$.
The map has derivative
\begin{alignat}{9} \frac{\partial\dens}{\partial\jj}_{|\spm,\Hpm} &= - \Big(\frac{\partial\Hpm}{\partial\dens}_{|\spm,\jj}\Big)^{-1}\frac{\partial\Hpm}{\partial\jj}_{|\spm,\dens}
&= -\frac{\jj^T}{\dens\csnd^2(1-\Mach^2)} .  \myeqlabel{eq:jj-dens-derivative}\end{alignat} 
For \emph{supersonic} flows \myeqref{eq:dHpm-ddens-atjj} is negative, so again we can obtain an analytic map $\jj\mapsto\dens$,
but the $\dens$ values are smaller, and strictly decreasing as $\vlen{\jj}$ decreases. 
By \myeqref{eq:ddens-dvv} $\dens$ is always decreasing as $\vlen\vv$ increases,
showing the well-known fact that for constant $\Hpm,\spm$ mass flux $\vlen\jj$ increases with speed $\vlen\vv$ in the subsonic range, but decreases in the supersonic range.

$\vpot$ solves \myeqref{eq:irrotationality} leaving \myeqref{eq:continuity} as the remaining scalar equation, for $\stf$ it is the other way around:
\begin{alignat}{9} 0 = \nabx^\perp\dotp(\dens^{-1}\jj) = -\nabx(\dens^{-1}\nabx\stf) \end{alignat}
is a 2nd order PDE for scalar $\stf$. 
Assuming smoothness, some calculation using \myeqref{eq:jj-dens-derivative} produces the non-divergence form
\begin{alignat}{9} 0 = \Aco:\nabx^2\stf ,\myeqlabel{eq:interior-stf}\end{alignat}
with same coefficient matrix as \myeqref{eq:vpot-interiorpde-xy}.

\mysection{Shocks and solid boundaries}

To consider shocks we use the standard shock relations for a smooth shock with downstream-pointing unit normal $\shockn$, with $\jmp{f}=f-f_i$ denoting the jump from upstream limit $f_i$ to downstream limit $f$ of some quantity $f$: 
\begin{alignat}{9} 0 &= \jmp{\jj}\dotp\shockn , \myeqlabel{eq:shock-jn}
\\ 0 &= \jmp{\vv}\dotp\shockn^\perp , \myeqlabel{eq:shock-ut}\end{alignat}
and implicitly $0 = \jmp{\Hpm} = \jmp{\spm}$, as we assumed them to be global constants. 
\defm{Admissible} shocks satisfy
\begin{alignat}{9} \densi < \dens ,\end{alignat}
which is equivalent to
\begin{alignat}{9} \csnd > \vv\dotp\shockn ,\end{alignat}
as well as
\begin{alignat}{9} \vvi\dotp\shockn > \csndi \end{alignat}
and
\begin{alignat}{9} \vvi\dotp\shockn > \vv\dotp\shockn. \end{alignat}

At solid walls with unit normal $\walln$ we use the usual \defm{slip condition}
\begin{alignat}{9} 0 = \jj\dotp\walln, \myeqlabel{eq:slip}\end{alignat}
which means $0=\walln^\perp\dotp\nabx\stf$ where $\walln^\perp$ is a tangent to the wall, so $\stf=\const$. Since we do not consider multiple disconnected walls there is no harm in taking
\begin{alignat}{9} 0 = \stf. \myeqlabel{eq:slip-stf} \end{alignat}

Using the stream functions $\stf,\stfi$ on the two sides of the shock, \myeqref{eq:shock-jn} implies that
\begin{alignat}{9} 0 = \shockn^\perp\dotp\jmp{\nabx\stf} = \shockn^\perp\dotp\nabx\jmp\stf ; \end{alignat}
the right-hand side is the tangential derivative of $\jmp\stf$, showing it is constant along the shock. In our case the shock meets a wall (or axis) where $\stf,\stfi$ satisfy the same slip condition \myeqref{eq:slip-stf}, showing the jump is zero along the shock:
\begin{alignat}{9} \stf = \stfi.  \myeqlabel{eq:shock-jmp-stf} \end{alignat}

For $\vpot$ the slip condition is
\begin{alignat}{9} 0 = \walln\dotp\nabx\vpot \end{alignat}
for wall normal $\walln$. The shock condition \myeqref{eq:shock-ut} $0=\jmp{\vv}\dotp\shockn^\perp=\shockn^\perp\dotp\nabx\jmp{\vpot}$ implies the tangential derivative of $\vpot$ is zero, so
\begin{alignat}{9} \vpot - \vpoti = \const \myeqlabel{eq:shock-jmp-vpot}\end{alignat}
along the shock. The constant need not be zero, and in this paper we do not need to fix its value.

\mysection{Shock polars}\mylabel{section:shockpolars}

We recall some prerequisites from \cite{elling-polarpotf} and other papers.
More general convex equations of state are considered there, but we only need to specialize to our $\gisen>1$ choice,
for which some of these properties have already been shown in \cite[Theorem 1]{elling-sonic-potf} and \cite[section 2.5]{elling-liu-pmeyer}.

\begin{figure}
\parbox{.49\textwidth}{%
  \centerline{\footnotesize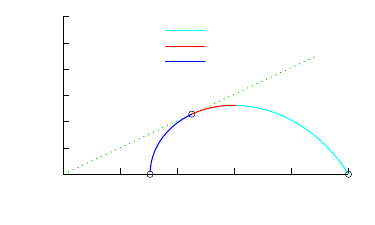}
  \caption{Velocity polar for $\Machi=2$, $\gisen=\frac53$.}
  \myfiglabel{fig:velocitypolar}
}\parbox{.49\textwidth}{%
  \centerline{\footnotesize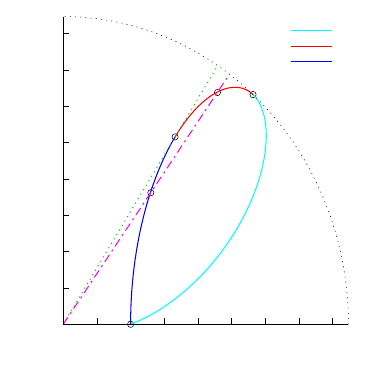}
  \caption{Mass flux polar for $\Machi=3$, $\gisen=\frac75$.}
  \myfiglabel{fig:massfluxpolar}
}
\end{figure}

We fix the pre-shock state $\vvi,\csndi,\densi$.

Let $\nuiavan=\arcsin(1/\Machi)$ be the upstream Mach angle. 
Let $\nuia\in\cli{-\nuiavan}{\nuiavan}$ be the angle from downstream-pointing unit shock normal $\shockn$ to $\vvi$, positive if counterclockwise. 
Angles outside that range produce non-admissible shocks.
For each $\nuia$ there is a unique solution $\vv$ of the shock relations \myeqref{eq:shock-jn} and \myeqref{eq:shock-ut};
$\nuia\mapsto\vv$ is analytic on $\cli{-\nuiavan}{\nuiavan}$, with $\vv=\vvi$ when $\nuia=\pm\nuiavan$ (vanishing point).
This curve is called \defm{$\vv$-shock polar}\footnote{In other texts some other curves, for example pressure-angle, are also called ``shock polar''.} (see fig.\ \myref{fig:velocitypolar}).
It can be defined as the set of $\vv$ with
\begin{alignat}{9} 0 = (\vv-\vvi)\dotp(\dens\vv-\densi\vvi),  \myeqlabel{eq:vv-polar-eqn} \end{alignat}
(with $\dens$ a function of $\vv$), since the shock conditions \myeqref{eq:shock-ut} and \myeqref{eq:shock-jn} imply that $\vv-\vvi$ is a shock tangent, $\dens\vv-\dens\vvi$ a shock normal.

Our existence proof in the second half of the paper will need $\pq=\nabx\stf=\jj^\perp$,
because the slip condition in the form $0=\walln\dotp\nabx\vpot$ would lead to rather awkward technical complications that can be avoided by using the stream function $\stf$ instead. 
The $\pq$-shock polar is simply the
\defm{$\jj$-polar} (mass flux polar) $\nuia\mapsto\jj=\dens\vv$ (see fig.\ \myref{fig:massfluxpolar}) rotated 90 degree counterclockwise around the origin $\veczero$.
The $\jj$-polar was already introduced in \cite{elling-massmodel} for other purposes; the needed properties follow trivially from those of the $\vv$-polar, so no separate proof is needed.
By \myeqref{eq:vv-polar-eqn} the $\jj$-polar can be defined as the set of $\jj$ with
\begin{alignat}{9} 0 = (\dens^{-1}\jj-\dens^{-1}\jji)\dotp(\jj-\jji),  \myeqlabel{eq:jj-polar-eqn} \end{alignat}
where $\dens$ is a function of $\jj$. 

For all curves $\nuia$ is a non-degenerate parameter, meaning $d\vv/d\nuia\neq 0$ and $d\jj/d\nuia\neq 0$.

By \myeqref{eq:shock-ut} $\vvi-\vv$ is normal to the shock, so $\nuia$ is also the counterclockwise angle from $\vvi-\vv$ to $\vv$.
For non-vanishing admissible shocks we have the downstream-pointing unit shock normal
\begin{alignat}{9} \shockn = \frac{\vvi-\vv}{\vlen{\vvi-\vv}}. \myeqlabel{eq:shockn-downstream}\end{alignat}
Exactly $\nuia=0$ produces a normal shock $\vv=\vvn\parallel\shockn\parallel\vvi$.
The $\jj$ polar is a bit different in that normal and vanishing point coincide, i.e.\ $\dens\vvn=\dens\vvi$ (see fig.\ \myref{fig:massfluxpolar}),
because \myeqref{eq:shock-jn} shows their normal parts are equal, whereas both tangential parts are zero.

$\vv$- and $\jj$-polar are symmetric under reflection across the line through $\veczero$ and $\vvi$. 
The \defm{upper halfpolar} is the $\nuia\geq 0$ part, its mirror image $\nuia\leq 0$ is the \defm{lower halfpolar}.
Results below for the upper halfpolar have obvious mirror images on the lower one.

There is a \emph{unique} $\nuiacrit\in\boi{0}{\nuiavan}$ so that the \defm{turning angle} $\turnap$ counterclockwise from $\vvi$ to $\vv$ (or $\jji$ to $\jj$)
is locally uniformly (l.u.) increasing\footnote{%
``Locally uniformly increasing'' can be defined for nonsmooth functions of course, but here simply refers to continuous derivatives along the curve being positive in some points.
} for $\nuia$ in the \defm{strong-type} range $\roi{0}{\nuiacrit}$, from $0$ to the \defm{critical angle} $\turnamax$,
but l.u.\ \emph{decreasing} for $\nuia$ in the \defm{weak-type} range $\boi{\nuiacrit}{\nuiavan}$, back to $0$.
The value $\vv$ (or $\jj$) for $\nuia=\nuiacrit$ is the \defm{critical point} $\vvc$ (or $\jjc$) of the shock polars. 

There is a \emph{unique} $\nuiasonic$ in the weak-type range $\boi{\nuiacrit}{\nuiavan}$ so that $\Mach=1$.
This \defm{sonic point} $\vv_1$ (or $\jj_1$) splits $\roi{0}{\nuiavan}$ into a supersonic range $\boi{\nuiasonic}{\nuiavan}$ and a subsonic range $\roi{0}{\nuiasonic}$ (elsewhere: ``transonic'').
The strong-type and critical shocks are always subsonic.
While velocity $\vlen{\vv}$ is l.u.\ increasing from normal to vanishing along the entire upper halfpolar, 
mass flux $\vlen\jj$ is l.u.\ increasing in the subsonic range, but l.u.\ decreasing in the supersonic range, with maximum $\jmax$ at the sonic point,
so on the $\jj$-polar the sonic point is especially obvious: the point farthest from $\jji$ (fig. \myref{fig:massfluxpolar}). 

\cite{elling-convexpolar} proves that $\vv$-polars for \emph{some} equations of state, including the $\gisen>1$ ones \myeqref{eq:polytropic-pressure},
are convex curves (see also \cite[Theorem 1]{elling-sonic-potf}),
a property that is comforting but not universal and not used in this paper. 
Useful results about shock polars for full Euler and other models can be found in \cite{teshukov-polar,elling-idealpolar,henderson-menikoff,elling-waalspolar,menikoff-plohr,elling-isenpolar}.

\mysection{Main theorems}\mylabel{section:mainresults}

\begin{theorem}
  \mylabel{th:weak-inf-univalent}%
  Given $\gisen\in\boi1\infty$, and constant $\densi>0$, $\csndi>0$ and $\Machi>1$, let $\vvi=(\Machi\csndi,0)$ and
  \begin{alignat}{9} \stfi \defeq \densi\vxi\ye \myeqlabel{eq:stfi-def}\end{alignat}
  the pre-shock stream function. 
  We consider $\hw\in\Cinf(\RR)$
  parametrizing\footnote{Since the wall is convex, and converging to a sub-critical slope at infinity, there is no limitation in parametrizing it by $\ye$.} the wedge wall
  \begin{alignat}{9}  \eWedgexy  \defeq \selset{ (\xe,\ye) }{ \ye\in\RR,~\xe=\hw(\ye) }  \end{alignat}
  with $\hw(0)=0$ (passing through the origin), mirror symmetry
  \begin{alignat}{9}  \hw(-\ye)=\hw(\ye)  \quad\text{ for $\ye\in\RR$,} \myeqlabel{eq:mirror-symmetry} \end{alignat}
  tangent converging at infinity like
  \begin{alignat}{9}  \frac{d\hw}{d\ye}(\ye) \conv \pm\cot\turnai \quad\text{as $\ye\conv\pm\infty$} \myeqlabel{eq:hw-slope}  \end{alignat}
  for some constant\footnote{%
    We only assume the tangent converges, which is strictly weaker than curve converging to a fixed line.}
  $\turnai\in\boi{0}{\turnamax}$, and locally uniformly convex wedge:
  \begin{alignat}{9} \frac{d^2\hw}{d\ye^2}(\ye) > 0 \quad\text{for $\ye\in\RR$.} \myeqlabel{eq:wedge-convex}\end{alignat}
  Let $\hs\in\Cinf(\RR)$ parametrize\footnote{There is no limitation in parametrizing the shock by $\ye$,
  since admissible shocks with $\vvi=(*,0)$ have tangent bounded away from horizontal.} the shock $ \eShockxy \defeq \selset{ (\xe,\ye) }{ \ye\in\RR,~\hs(\ye)=\xe }$ with
  \begin{alignat}{9}  \hs(\ye) < \hw(\ye)  \quad\text{for $\ye\in\RR$}, \end{alignat}
  and define the post-shock region $\eDomxy \defeq \selset{ (\xe,\ye) }{ \ye\in\RR,~\hs(\ye)<\xe<\hw(\ye) }$.
  Let $\vpot,\stf\in\Cinf(\ecDomxy)$ represent a subsonic irrotational Euler flow solving \myeqref{eq:interior-stf} in $\eDomxy$,
  mirror-symmetric in the sense $\stf(\xe,\ye)=-\stf(\xe,-\ye)$ for all $\xx\in\eDomxy$,
  satisfying
  the slip condition \myeqref{eq:slip-stf} on $\eWedgexy$
  and the shock relations \myeqref{eq:shock-jmp-stf} and \myeqref{eq:shock-jmp-vpot} across $\eShockxy$.
  We call such a flow \defm{weak-at-infinity} if
  \begin{alignat}{9}  \jj\conv\jjw \quad\text{as $\ye\conv+\infty$ while $\xx\in\ecDomxy$,} \myeqlabel{eq:weak-at-inf}\end{alignat}
  strong-at-infinity for $\jj\conv\jjs$.

  Under these assumptions, any weak-at-infinity flow is \defm{univalent}, meaning $\ecDomxy\ni\xx\mapsto\nabx\stf(\xx)$ is an injective map.
\end{theorem}

\begin{remark}
  \mylabel{rem:analyticity}%
  We could assume less regularity than $\Cinf$,
  but standard Schauder estimates improve regularity in the interior $\eDomxy$ to $\Cinf$ anyway (e.g.\ \cite[Theorem 6.17]{gilbarg-trudinger});
  in fact as the coefficients $\Aco=\idm-(\vv/\csnd)^2$ of our interior PDE \myeqref{eq:interior-stf} $0=\Aco:\nabx^2\stf$ are analytic,
  in particular $\csnd$ \myeqref{eq:csnd} via $\dens$ \myeqref{eq:bernj},
  the solution $\stf$ is also analytic there, by classical results of elliptic PDE theory, see for example the general treatment in \cite[section 6.7]{morrey-book}
  (also \cite{bernstein-1904,eberhard-hopf-elliptic-analytic-1932,morrey-analyticity-i,friedman-analyticity-jmathmech1958} for historical perspective).
  
  The shock relations are also defined by analytic functions, so the shock and $\stf$ at the shock are analytic.
  Although this has been known for a long time and used in classical papers, it is apparently not universally known,
  so in the appendix section \myref{section:shock-analyticity} we supply a short proof by reduction from free to fixed boundary,
  which permits using standard results on analyticity (e.g.\ \cite[section 6.7]{morrey-book}; all our problems satisfy proper ellipticity and the complementing condition),
  but also Schauder estimates for bootstrapping shock regularity from $\Ck{m,\alpha}$ to $\Cinf$.

  For the wedge $\eWedgexy$ we consciously do \emph{not} require analyticity, which is a controversial assumption due to its non-local nature.
\end{remark}

Theorem \myref{th:weak-inf-univalent} shows that every weak-at-infinity subsonic post-shock is univalent. 
It does \emph{not} show such a flow must exist; this cannot be expected in general,
as it is a well-known phenomenon of aerodynamics that subsonic flight permits supersonic regions to appear at the airfoil where the surface is strictly convex.
This phenomenon is easily observed for all kinds of model equations including irrotational Euler, 
not only for airfoils but also for convex blunt wedges and all kinds of other applications.
Mathematically it is reflected in the failure of maximum principles for Mach number at strictly convex slip-condition boundaries. 
  
However, Theorem \myref{th:weak-inf-existence} proves that for any $\gisen,\Machi$ there exists an infinity of smooth blunt convex wedge shapes
allowing flows with weak-at-infinity shocks with subsonic post-shock region.

For strong-at-infinity flows there are two distinct cases, of which one is not univalent (see Remark \myref{rem:strong-at-inf});
this requires a complicated discussion left for future work.

We state the following existence result as flow along a wall with convex ramp (fig.\ \myref{fig:psilevel}), which is one among several closely related problems of independent interest;
existence of flow onto symmetric convex blunt wedges follows by a simple odd reflection to the lower halfplane:
\begin{alignat}{9} \stf(\xe,-\ye)\defeq-\stf(\xe,\ye)  \quad\text{for $\ye\geq 0$, $\xe\in\RR$.} \end{alignat}

\begin{theorem}
  \mylabel{th:weak-inf-existence}%
  For any $\gisen\in\boi1\infty$, $\densi>0$, $\csndi>0$, $\Machi>1$ and $\turnai\in\boi{0}{\turnamax}$
  there are infinitely many $\hw\in\Conce\roi{0}{\infty}\isect\Cinf\boi{0}{\infty}$ so that
  \begin{alignat}{9} \hw(0) = 0 \ ,\ \frac{d\hw}{d\ye}(0) = 0 \ ,\ \frac{d^2\hw}{d\ye^2}(\ye) > 0 \ ,\  \frac{d\hw}{d\ye}(\ye) \conv \arccot\turnai \quad(\ye\in\RR), \end{alignat}
  together with $\hs\in\Conce\roi{0}{\infty}\isect\Cinf\boi{0}{\infty}$ so that
  \begin{alignat}{9} \hs(\ye) < \hw(\ye) \quad\text{for $\ye\in\roi{0}{\infty}$,} \end{alignat}
  and with subsonic $\stf\in\Conce(\cDomxy)\isect\Cinf(\cDomxy\setdiff\sett{\xxn,\xxo})$ 
  for normal point $\xxn=(\hs(0),0)$ and post-shock domain $\Domxy=\selset{(\xe,\ye)}{\ye>0,\hs(\ye)<\xe<\hw(\ye)}$,
  solving the irrotational Euler equations \myeqref{eq:interior-stf} in $\Domxy$,
  the slip condition \myeqref{eq:slip-stf} on $\Wedgexy=\selset{(\xe,\ye)}{\ye>0,\xe=\hw(\ye)}$ and $\Axisxy\defeq\selset{(\xe,\ye)}{\ye=0,\hs(\ye)<\xe<\hw(\ye)}$
  and the shock conditions \myeqref{eq:shock-jn} and \myeqref{eq:shock-ut} across $\Shockxy=\selset{(\xe,\ye)}{ \hs(\ye)=\xe} $.
  The flow is univalent, meaning $\cDomxy\ni\xx\mapsto\nabla\stf(\xx)$ is an injective map, and it is weak-at-infinity in the sense
  \begin{alignat}{9} \jj \conv \jjw \quad\text{as $\vlen{\xx}\conv\infty$ with $\xx\in\Domxy$.} \end{alignat} 
\end{theorem}
\begin{remark}
  In particular we can construct examples for $\Machi-1>0$ arbitrarily small,
  so that any admissible shock has small strength everywhere;
  small irrotational Euler shocks are very close to full Euler shocks.
  The post-shock region still has large variation of velocities ($\vv$ almost sonic in the normal point $\xxn$, zero at the stagnation point $\xxo$).
  Therefore this paper has not only the first example of weak-at-infinity blunt-wedge flows, but also the first example with arbitrarily small shocks.
\end{remark}

For theorem \myref{th:weak-inf-univalent} the main obstacle is that because the wedge is strictly convex the fluid domain is \emph{non-convex}.
It is a well-known obstacle in elliptic PDE that non-convex domains have far fewer maximum principles than convex ones. 
A famous classical case are gradients for the minimal-surface equation.
It can be considered the $\gisen\dnconv-1$ limit of the irrotational Euler equations,
for which Mach numbers satisfy well-known maximum principles at boundary points where the domain is locally convex, but not otherwise. 

In our case we are forced to work with stream function and velocity angle. The angle cannot even be defined until we know stagnation points do not occur away from the wedge,
which is proven in section \myref{section:streamfunctionlevels} by analyzing stream function level sets in all detail.
The angle has rather complicated level sets with some singular behaviour and branching,
but after long analysis in section \myref{section:vanglevels} enough can be shown to prove that weak-at-infinity shocks are necessarily univalent.
(This is wrong for strong-type infinity shocks, demonstrating that the result is delicate
and that the importance of types remains neglected and misunderstood.)

Every flow that exists is necessarily univalent, so for theorem \myref{th:weak-inf-existence} we might as well work in the $\vv$-plane. 
However, this would cause some hidden technical obstacles that may be discussed in future work. 
Instead we find that an approach can be made by working in the unusual mass flux $\jj$-plane.

The existence requires solving an oblique derivative problem with some troublesome properties,
including oblique degeneracy and leading-order parts of operators vanishing in some points.
We solve the problem for all parameters, so there is no smallness to use;
proving existence requires a detailed estimation of corner exponents and asymptotics in section \myref{section:corner-exponents}
followed by a long analysis of Fredholm indices for operators on weighted Sobolev spaces in section \myref{section:existence-pq}.
Given the $\jj$-plane solution it is necessary to return to original coordinates, in section \myref{section:pq-to-xx}.
Again lack of smallness forces us to develop some new techniques to prove injectivity for the coordinate change,
in particular certain maximum principles in the $\jj$ plane, again in domains that need not be convex.
We note that even the $\vv$-shock polar need not be convex for some equations of state \cite{elling-convexpolar}.

\mysection{General level set result}

We recall and extend some technical results, some of which were already used in \cite[Lemma 3.1--3.3]{elling-liu-pmeyer} as well as \cite{elling-polylow,elling-threecorner-arxiv,elling-threecorner} for rather different purposes.

Recall that every non-singular analytic curve has
an analytic unit-speed parametrization $s\mapsto\vec f(s)$, i.e.\ with $\vlen{\vec f_s}=1$.

\newcommand{\genAsym}{M}
\newcommand{\genA}{\mat{\genAsym}}
\newcommand{\genvav}{m}
\newcommand{\genva}{\vec\genvav}
\newcommand{\genaz}{\genvav^0}
\newcommand{\genbv}{b}
\newcommand{\genbb}{\vec\genbv}
\newcommand{\genbzero}{\genbv^0}

\begin{proposition}
  \myeqlabel{prop:levels-atpoint}%
  
  In some Cartesian coordinates $(\tau,\nu)$ consider
  \begin{alignat}{9} \genA:\nabla^2\genphi + \genva\dotp\nabla\genphi = 0 , \myeqlabel{eq:genell}\end{alignat}
  where $\genphi,\genva,\genA$ are real analytic functions of $(\tau,\nu)$, with values scalars, vectors, symmetric positive definite matrices. 
  Assume there is\footnote{else $\genphi$ is locally constant} a least positive integer $m$ with $D^m\genphi(0,0)\neq 0$. 
 
  (i) If \myeqref{eq:genell} holds on the disk $V$ of radius $\delta>0$ around $(0,0)$, then for $\delta>0$ sufficiently smaller the level set
  \begin{alignat}{9} L=\selset{(\tau,\nu)\in V}{\genphi(\tau,\nu)=\genphi(0,0)} \end{alignat}
  is a union of $2m$ analytic curves with one endpoint $(0,0)$ where their inward unit tangents are pairwise distinct.

  (ii) Consider the boundary $V_0$ defined as the graph of an analytic function $h:\boi{-\delta^2}{\delta^2}\rightarrow \boi{-\delta}{\delta}$ for $\delta>0$:
  for $\sigma\in\sett{-,0,+}$, 
  \begin{alignat}{9} V_\sigma := \selset{(\tau,\nu) \in  \boi{-\delta^2}{\delta^2}\cartp\boi{-\delta}{\delta} }{ \sign(\nu-h(\tau))=\sigma } . \end{alignat}
  Assume $h(0)=0$ and $h'(0)=0$ (so $\partial_\tau$ is tangential).
  If \myeqref{eq:genell} holds on $V=V_0\union V_+$, 
  then for $\delta>0$ sufficiently smaller
  \begin{alignat}{9} L=\selset{(\tau,\nu)\in V_0\union V_+}{\genphi(\tau,\nu)=\genphi(0,0)} \end{alignat}
  is a union of $m-1$ or $m$ or $m+1$ analytic curves with endpoint $(0,0)$ where their inward unit tangents are pairwise distinct. 

  (iii) In (ii) assume additionally the boundary condition
  \begin{alignat}{9} 0 = (\dbb\dotp\nabla+\bzero)\genphi \quad\text{on $V_0$} \myeqlabel{eq:gen-bcond}\end{alignat}
  for analytic $\dbb$ and \emph{constant} $\bzero$. If $\dbb$ is not tangential,
  then in (ii) the union has exactly $m$ curves, with inward unit tangents in $(0,0)$ not tangential to $V_0$. 

  (iv) Across each curve the sign of $\genphi-\genphi(0,0)$ changes to the opposite. 
\end{proposition}
\begin{proof}
  For convenience we first change to new coordinates $(\xi,\eta)$ so that $\genA(\veczero)=\idm$ %(
  (the Jacobian of the coordinate transformation is an inverse square root of the symmetric positive definite $\genA(\veczero)$.)
  Consider $m\geq 2$; for $j=0,...,m-2$ take $T=\partial_\xi^{m-2-j}\partial_\eta^j$ of \myeqref{eq:genell}:
  \begin{alignat}{9} 0 = (\genA:\nabla^2)T\genphi + [T,\genA:\nabla^2]\genphi + (\genva\dotp\nabla)T\genphi + [T,\genva\dotp\nabla]\genphi . \myeqlabel{eq:Acomm}\end{alignat}
  $(\genva\dotp\nabla)T$ and the commutators $[T,\genA:\nabla^2]\genphi,[T,\genva\dotp\nabla]$ are operators of order $\leq m-1$,
  so applied to $\genphi$ they yield zero in the origin where $D^k\genphi(\veczero)=0$ for $1\leq k\leq m-1$, so that using $\genA=\idm$ there we find
  \begin{alignat}{9} 0 &= (\genA:\nabla^2)T\genphi = (\partial_\xi^2+\partial_\eta^2)\partial_\xi^{m-2-j}\partial_\eta^j\genphi   \qquad(j=0,...,m-2)
  \notag\\\impl\quad
  \alpha &:= \partial_\xi^m\genphi = - \partial_\xi^{m-2}\partial_\eta^2\genphi = + \partial_\xi^{m-4}\partial_\eta^4\genphi = ... ,
  \notag\\ -\beta &:= \partial_\xi^{m-1}\partial_\eta\genphi = - \partial_\xi^{m-3}\partial_\eta^3\genphi = + \partial_\xi^{m-5}\partial_\eta^5\genphi = ... \quad, \notag\end{alignat}
  which we can write in compact form as
  \begin{alignat}{9} \partial_\xi^{m-j}\partial_\eta^j\genphi(\veczero) = \Re(\subeq{(\alpha+i\beta)}{\eqdef c}i^j) \qquad(j=0,...,m). \end{alignat}
  By definition of $m$ we have $c\neq 0$. 
  The degree $m$ Taylor polynomial $T_m\genphi$ at $\veczero$ with remainder term $R_m\genphi$ is
  \begin{alignat}{9} T_m\genphi
  &= \sum_{j=0}^m \binom{m}{j} \partial_\xi^{m-j}\partial_\eta^j\genphi(\veczero) \xi^{m-j} \eta^j
  \\&= \Re\big( c \sum_{j=0}^m \binom{m}{j} i^j \xi^{m-j} \eta^j \big)
  \\&= \Re\big( c (\subeq{\xi+i\eta}{\eqdef z})^m \big). \myeqlabel{eq:Tm}\end{alignat}
  Using polar forms of $c,z$ quickly shows
  the zero set is a union of $m$ lines through $(\xi,\eta)=(0,0)$ at equal angles $\pi /m$.
  The condition $\Aco(\veczero)=\idm$ for choosing coordinates $(\xi,\eta)$ allows an arbitrary rotation around the origin,
  so we may assume none of the lines is vertical,
  so that
  \begin{alignat}{9} T_m\genphi = \xi^m \subeq{\Re(c(1+is)^m)}{\eqdef P_m(s)}  \quad\text{with}\quad  s \defeq \eta/\xi  \end{alignat}
  is a degree $m$ polynomial in $s$, with \emph{finite} distinct real roots $s_1,...,s_m$.

  Since $\genphi$ is analytic we get
  \begin{alignat}{9} \genphi = 
  \xi^m ( \subeq{ P_m(s) + \xi U_m(s,\xi) }{=:h_m(s,\xi) } ) \end{alignat}
  for some analytic function $U_m$ of $\xi,s$ extracted from the remainder $R_m$ by using $\xi^{m+j-k}\eta^k=\xi^m\cdot \xi\cdot \xi^{j-1}s^k$ for each monomial of degree $m+j$ with $j\geq 1$. 
  Clearly
  \begin{alignat}{9} \partial_sh_m(s_\ell ,0) = \partial_sP_m(s_\ell ) + 0\cdot \partial_sU_m(s_\ell ,0) \neq  0 \end{alignat}
  since the roots $s_1,...,s_m$ are distinct. 
  So the implicit function theorem for analytic functions \cite[thm.\ 2.3.1]{krantz-parks} gives for each $\ell $ an analytic $\xi\mapsto S_\ell (\xi)$ with $S_\ell (0)=s_\ell $ solving
  \begin{alignat}{9} h_m(S_\ell (\xi),\xi) = 0. \end{alignat}
  The graphs of $\xi\mapsto \eta=\xi S_\ell (\xi)$ split into $\xi\leq 0$ and $\xi\geq 0$ halves that give, after returning to original coordinates, the $2m$ curves in the statement. 

  (ii) 
  $\genA,\genva,\genphi$ are analytic on $V_0\union V_+$, so they have analytic extensions to $V_-$ (for $\delta>0$ sufficiently small).
  Since $\genA:\nabla^2\genphi+\genva\dotp\nabla\genphi$ is zero on $V_0\union V_+$, by analyticity it is still zero on $V_-$.
  So we can repeat the proof of (i), obtaining $m$ analytic curves \emph{through} $(\tau,\nu)=(0,0)$, 
  with pairwise distinct tangents there. If one of the curves is tangent to the boundary $V_0$, 
  let $\nu=g(\tau)$ parametrize it, then by analyticity $h-g$ has isolated zeros, so for $\tau\in\boi{0}{\delta}$ with $\delta>0$ sufficiently smaller $\tau\mapsto \sign(g-h)$ is a constant $\sigma$, with the right half-curve entirely contained in $V_\sigma$; for another constant $\sigma$ the same holds for the left half-curve. The $m-1$ other curves, or all $m$ curves if none tangent to $V_0$, have one half entirely in $V_-$, the other half in $V_+$, again for $\delta>0$ sufficiently smaller.

  (iii) 
  Since the change to $(\xi,\eta)$ coordinates allows an arbitrary rotation around $(0,0)$, we align $\partial_\xi$ to be tangential to $V_0$.
  Take $\partial_\xi^{m-1}$ of \myeqref{eq:gen-bcond} and evaluate in $(\xi,\eta)=(0,0)$, using $\nabla^k\genphi=0$ there for $k=1,...,m-1$ (recall $\bzero$ constant):
  \begin{alignat}{9} 0
  &= \bbsym^\xi\partial_\xi^m\genphi+\bbsym^\eta\partial_\xi^{m-1}\partial_\eta\genphi
  \\&= \bbsym^\xi\alpha-\bbsym^\eta\beta.
  \end{alignat}
  By assumption $\dbb$ is not tangential, i.e.\ $\bbsym^\eta\neq 0$, so $\alpha+i\beta\neq 0$ implies $\alpha\neq 0$. 
  This means $T_m\genphi\topref{eq:Tm}{=}\Re((\alpha+i\beta)(\xi+i\eta)^m)$ does not have zeros $(\xi,\eta)$ with $\eta=0$,
  so none of the $m$ zero lines is horizontal, and since they are tangent to the $m$ curves through $(0,0)$ from (ii),
  none of those curves is tangential to $V_0$ and therefore one half of each passes into $V_+$. 

  (iv) This follows from the degree $m$ polynomial $P_m$ above having  $m$ distinct roots, necessarily with a sign change across each. 
\end{proof}

\begin{definition}
  Given a scalar analytic function $\genphi:V\rightarrow \RR$ on an open set $V\subset\RR^2$,
  in our case $V=\Domxy$ or $V=\Dompq$. 
  A point $\xx\in V$ is called \defm{vertex} if $\nabx\genphi(\xx)=0$.
  If it is not a vertex, then the level set through $\xx$ is locally a single non-singular analytic curve.
  We call such a curve an \emph{edge} if it is maximal (under inclusion) among analytic curves in $V$ through $\xx$ whose interior
  (excluding endpoints) does not contain vertices.
\end{definition}

\begin{remark}\mylabel{rem:single-endpoint-boundary}
  Let a scalar nonconstant analytic $\genphi:V\rightarrow \RR$ solve a PDE of type \myeqref{eq:genell} in each point of an open set $V$. 
  We will not need or discuss that some of the following properties also hold for arbitrary analytic or even $\Conce$ functions. 

  Due to $\nabx\genphi\neq 0$ the interior of an edge cannot intersect itself or interiors of other edges.
  The vertices are isolated, 
  and the edges in any compact subset $K$ of $V$ have finite total arclength (cover $K$ with the neighbourhoods from Proposition \myref{prop:levels-atpoint}(i) and pick a finite subcover).
  An end of an edge that does not converge to a vertex must leave any $K$ (converging either to $\bdry V$ or to infinity).

  Each level set can be thought of as an embedded planar graph (although the strict definition may require omitting edges that do not converge to boundary or infinity at any end).
  
  On connected \emph{analytic} segments of the boundary, each level set is either the entire segment or a discrete subset. 
  If some end of an edge approaches the segment, then by Proposition \myref{prop:levels-atpoint}(ii) it must converge to exactly one point of the segment.

  At non-analytic boundary segments, or in corners, 
  we encounter only cases where the level set is known to be well-behaved based on additional information. 
\end{remark}

\begin{proposition}
  \mylabel{prop:levels-tree}%
  If a non-constant analytic $\genphi:V\rightarrow \RR$ solves \myeqref{eq:genell}, then the connected components of each level sets are ``trees'', meaning they do not contain piecewise analytic simple loops. 
\end{proposition}\begin{proof}
  If some level set contains a piecewise analytic simple loop, then that encloses by the Jordan curve theorem a bounded open nonempty set with compact closure $K$.
  By compactness and continuity $\genphi$ attains minima and maxima in $K$, but the strong maximum principle for \myeqref{eq:genell} does not permit interior extrema for non-constant $\genphi$, and if both minimum and maximum occur on $\bdry K$, which is part of a level set, then $\genphi$ is constant on $K$ which has nonempty interior, again a contradiction. 
\end{proof}

To decide whether a point is a vertex the following simple lemma is useful.
\begin{lemma}
  \mylabel{lemma:two-zero-three}%
  Let $\genphi$ solve $ 0 = \Aco:\nabla^2\genphi $ (\emph{without} first-order term) in a point where\\
  (i) $\dex_1\dotp\nabla^2\genphi\,\dex_2=0=\dex_1\dotp\nabla^2\genphi\,\dex_2'$ for $\dex_1$ nonzero, $\dex_2,\dex_2'$ linearly independent, or \\
  (ii) $\det\nabla^2\genphi=0$. \\
  Then $\nabla^2\genphi=0$ there. 
\end{lemma}\begin{proof}
  (i) means $\dex_1\dotp\nabla^2\genphi=0$, but then $\nabla^2\genphi$ has an eigenvalue zero, implying (ii), which implies $\nabla^2\genphi$ is a rank $\leq 1$ matrix $\vec u^2$ (symmetric),
  so $0=\Aco:\nabla^2\genphi=\vec u\dotp\Aco \vec u$, but $\Aco$ positive definite so $\vec u=0$. 
\end{proof}
Of course this is one of many properties of \emph{planar} elliptic PDE widely used in the classical literature
on quasi-conformal maps and various generalizations of analytic functions. 
Using those tools some of our steps and results could be recast in the language of Riemann surfaces,
but the relevant background is a bit scattered and published in irregular form \cite{bers-lecture-notes},
so it is easier and more readable to write self-contained elementary proofs.

\mysection{Stream function level sets}\mylabel{section:streamfunctionlevels}

We begin the proof of Theorem \myref{th:weak-inf-univalent}. 
We focus on the upper halfplane subsets $\Domxy,\Shockxy,\Wedgexy,\Axisxy$ as in Theorem \myref{th:weak-inf-existence};
results in the lower halfplane follow by mirror symmetry.

\begin{figure}
  \centerline{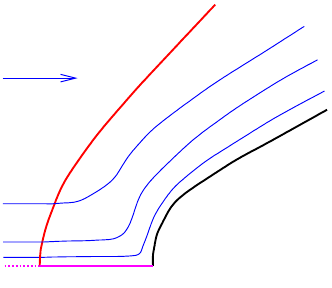}
  \caption{Streamfunction $\stf$ level sets are simple curves}
  \myfiglabel{fig:psilevel}
\end{figure}

\begin{proposition}
  \mylabel{prop:stf-levels-at-shock}%
  Each streamline (level set of $\stf$) meets the shock exactly once, is a single analytic curve near there, and uniformly separated from the shock everywhere else.
\end{proposition}\begin{proof}
  At $\Shockxy$ we have $\stf\topref{eq:shock-jmp-stf}{=}\stfi\topref{eq:stfi-def}{=}\densi\vxi\ye$ which is uniformly increasing in $\ye$,
  so each positive $\stf_+$ occurs exactly once;
  besides $\nabx\stf=\dens\vv^\perp\neq 0$ at the downstream side of any shock, so by the analytic implicit function theorem the level set is locally an analytic curve. 
\end{proof}

Near upper post-shock infinity each level set of $\stf$ is a single analytic curve running to infinity; more precisely:
\begin{proposition}
  \mylabel{prop:stf-levels-near-inf}%
  For any $\stf_+\in\boi0\infty$ there is a $\ye_+\geq 0$ and an analytic function  $h:\boi{\ye_+}{\infty}\rightarrow\RR$
  so that for $(\xe,\ye)\in\cDomxy$ with $\ye>\ye_+$,
  \begin{alignat}{9} \stf(\xx)=\stf_+ \qeq \xe=h(\ye). \end{alignat} 
\end{proposition}\begin{proof}
  By assumption $\stf_\xe=-\jy\topref{eq:weak-at-inf}{\conv}-\jyw<0$ as $\ye\conv\infty$ in $\Domxy$,
  so $\stf$ is strictly monotone along any sufficiently high horizontal line.
  On the other hand, at the shock
  \begin{alignat}{9} \stf\topref{eq:shock-jmp-stf}{=}\stfi\topref{eq:stfi-def}{=}\densi\vxi\ye\upconv+\infty   \quad\text{as $\ye\upconv+\infty$,} \end{alignat}
  while $\stf=0$ at the wedge, so if the line is sufficiently high there is exactly one point with $\stf=\stf_+$. 
  Since $\stf$ is analytic with $\nabla\stf\neq 0$, the implicit function theorem for analytic functions shows existence of $h$.
\end{proof}

\begin{proposition}\mylabel{prop:stf-positive-level-sets}%
  $\Domxy$ is partitioned into level sets of $\stf$ for positive values $\stf_+$, each of which is a simple analytic curve from $\Shockxy$ to upper downstream infinity, and 
  \begin{alignat}{9} \nabx\stf &\neq 0 \qquad\text{on $\cDomxy\setdiff\sett{\xxo}$,} \myeqlabel{eq:nabx-stf-nowhere-zero} \end{alignat}
  and with $\partial_n$ normal derivative \emph{inward} to $\Domxy$
  \begin{alignat}{9} \stf_n &> 0 \quad\text{on $\Wedgexy\union\Axisxy$.} \myeqlabel{eq:normal-deriv-stf-walls}\end{alignat} 
\end{proposition}\begin{proof}
  Fix $\stf_+$. Its level set is a disjoint union of locally finite trees,
  by Proposition \myref{prop:levels-tree} applied to the elliptic interior PDE \myeqref{eq:interior-stf}.
  Any tree that has a vertex has at least four simple paths from that vertex, pairwise disjoint except at the endpoints, to boundary or to infinity.
  But they cannot approach $\Wedgexy\union\Axisxy$ because at those boundaries the slip condition $\stf=0$ requires separation from $\sett{\stf=\stf_+}$ for $\stf_+>0$, 
  whereas by Propositions \myref{prop:stf-levels-near-inf} and \myref{prop:stf-levels-at-shock} 
  there is exactly one edge to $\Shockxy$ and exactly one edge to upper post-shock infinity,
  not enough for more than two simple paths. 
  This contradiction shows there are no vertices, and by the same argument there is at most one tree, so the level set is a single edge running from $\Shockxy$ to infinity. 

  If $\nabx\stf(\xx)$ was zero in some $\xx\in\Domxy$, then there would be a smallest $m\geq 2$ with $\nabx^m\stf(\xx)\neq 0$,
  so Proposition \myref{prop:levels-tree} would show $\xx$ is a vertex, which we have already ruled out.
  The shock relations do not permit $\dens\vv^\perp=\nabx\stf=0$ at the shock either.
  $\stf$ attains a minimum value $0$ over $\cDomxy$ on $\Wedgexy$ and $\Axisxy$, which are smooth boundaries, and $\stf$ solves a scalar elliptic 2nd order PDE,
  so the Hopf lemma implies the outward normal derivative of $\stf$ is negative everywhere
  on $\Wedgexy\union\Axisxy$. This completes the proof of \myeqref{eq:nabx-stf-nowhere-zero} and \myeqref{eq:normal-deriv-stf-walls}.
\end{proof}

\mysection{Velocity angle level sets}\mylabel{section:vanglevels}

By \myeqref{eq:nabx-stf-nowhere-zero} $\jj=-\nabx^\perp\stf\neq 0$ everywhere in $\ecDomxy\setdiff\sett{\xxo}$, so $\vv\neq 0$ everywhere, 
so we may consider the velocity angle $\vang$.
Angles are by nature multi-valued, defined up to multiples of $2\pi$, so care is needed to avoid subtle errors. 
The domain $\ecDomxy\setdiff\sett{\xxo}$ is simply connected, so we can well-define $\vang$ by requiring it to be continuous, satisfy
\begin{alignat}{9} \vv = \vlen{\vv} (\cos\vang,\sin\vang), \end{alignat} 
and fix
\begin{alignat}{9} \vang = 0  \quad\text{on $\Axisxy\union\sett{\xxn}$,} \myeqlabel{eq:vang-axis}\end{alignat}
where we have already shown $\vx=\dens^{-1}\stf_\ye>0$ \myeqref{eq:normal-deriv-stf-walls}, while $\vy=-\dens^{-1}\stf_\xe=0$ by the slip condition.

The discussion of level sets of velocity angle is much more complicated than that of $\stf$ level sets, because certain extrema \emph{can} occur along or at the shock,
and near $\xxo$ we do not have continuity.

Velocity angle satisfies a 2nd order elliptic PDE of its own, which is convenient because a lot of theory applies, for example the strong maximum principle or the Hopf lemma.

The Euler equations have symmetry under translation and rotation, so we may change (preserving orientation) to arbitrary Cartesian coordinates $(\nu,\tau)$ so that $\vpot_\nu$ is locally positive. (At the downstream side of shocks, with $\partial_\nu$ normal pointing downstream, $\vpot_\nu$ is always positive.)
Then
\begin{alignat}{9}  \vang  =  \const + \arctan\frac{\vpot_\tau}{\vpot_\nu}  ,  \myeqlabel{eq:vang-arctan}  \end{alignat}
where subscripts $\nu,\tau$ are partial derivatives, and the constant only depends on the choice of coordinates.
We calculate
\begin{alignat}{9} \nabla\vang &= \frac{ \vpot_\nu\nabla\vpot_\tau - \vpot_\tau\nabla \vpot_\nu }{ \vlen{\nabla\vpot}^2 }, \myeqlabel{eq:nabla-vang}
\\ \nabla^2\vang &= \frac{ \vpot_\nu\nabla^2\vpot_\tau - \vpot_\tau\nabla^2 \vpot_\nu }{ \vlen{\nabla\vpot}^2 } + Q_1 , \end{alignat}
where $Q_i=(\mat{Q}_i\nabla\vang+\vec{\ell }_i)\dotp \nabla\vang$ with $\mat{Q}_i,\vec{\ell }_i$ analytic functions of $\nabla\vpot$.
\begin{alignat}{9} \Aco(\nabla\vpot):\nabla^2\vang
&= \frac{\vpot_\nu\Aco(\nabla\vpot):\nabla^2\vpot_\tau-\vpot_\tau\Aco(\nabla\vpot):\nabla^2\vpot_\nu}{\vlen{\nabla\vpot}^2} + Q_2
\\&= \frac{\vpot_\nu(\supeq{\Aco(\nabla\vpot):\nabla^2\vpot}{=0})_\tau-\vpot_\tau(\supeq{\Aco(\nabla\vpot):\nabla^2\vpot}{=0})_\nu}{\vlen{\nabla\vpot}^2} + Q_3 , \myeqlabel{eq:ADDvang}\end{alignat}
because 
$(\Aco(\nabla\vpot))_k:\nabla^2\vpot$ for $k=\tau,\nu$ are also of $Q_i$ form. 
We combine \myeqref{eq:nabla-vang} and the $\vpot$ equation \myeqref{eq:vpot-interiorpde-xy} into 
\begin{alignat}{9} \begin{bmatrix} \Acosym^{\nu\nu}   & 2\Acosym^{\nu\tau}   & \Acosym^{\tau\tau}
  \\  -\vpot_\tau/\vlen{\nabla\vpot}^2 & \vpot_\nu/\vlen{\nabla\vpot}^2  & 0
  \\  0           & -\vpot_\tau/\vlen{\nabla\vpot}^2 & \vpot_\nu/\vlen{\nabla\vpot}^2  \end{bmatrix}
\begin{bmatrix}   \vpot_{\nu\nu} \\   \vpot_{\nu\tau} \\   \vpot_{\tau\tau} \end{bmatrix}
= \begin{bmatrix}   0 \\   \vang_\nu \\   \vang_\tau \\ \end{bmatrix} \myeqlabel{eq:vang-3x3}\end{alignat}
The system determinant is
\begin{alignat}{9} = \frac1{\vlen{\nabla\vpot}^4} \big( \Acosym^{\nu\nu} \vpot_\nu^2 + 2\Acosym^{\nu\tau} \vpot_\nu \vpot_\tau+ \Acosym^{\tau\tau} \vpot_\tau^2 \big)
= \frac1{\vlen{\nabla\vpot}^4} \nabla\vpot\dotp\Aco \nabla\vpot > 0
\myeqlabel{eq:vangsysdet}
\end{alignat}
for any nonzero $\nabla\vpot$ and any positive definite $\Aco$. 
So there is a unique solution $\nabla^2\vpot$, linear in $\nabla\vang$; we can substitute it into $Q_3$ of \myeqref{eq:ADDvang}, obtaining 
\begin{alignat}{9} \Aco:\nabla^2\vang + (\mat{Q}\nabla\vang+\vec{\ell })\dotp \nabla\vang , \myeqlabel{eq:vang-pde}\end{alignat}
still for $\mat{Q},\vec{\ell }$ analytic in $\nabla\vpot$;
of course for $\vang$ they have well-known explicit formulas, but their details are not needed.

\begin{proposition}\mylabel{prop:vang-interior-extrema}
    Velocity angle $\vang$ does not have local extrema in the interior $\Domxy$.
\end{proposition}\begin{proof}
  The strong maximum principle for \myeqref{eq:vang-pde} does not permit interior extrema unless $\vang=\measuredangle\vv$ is constant in the (path-connected) set $\Domxy$,
  which is clearly incompatible with the slip condition $\vv\dotp\walln=0$ at the wedge $\Wedgexy$ which is strictly convex \myeqref{eq:wedge-convex}.
\end{proof}

\begin{proposition}
  \mylabel{prop:angle-levels-stag}%
  For each $\vang_+\in\boi{0}{\pi/2}$ the level set $\selset{\xx\in\Domxy}{\vang(\xx)=\vang_+}$ is near $\xxo$ a single analytic curve,
  with tangent direction converging to angle $\pi-\vang_+$ as $\xx$ approaches $\xxo$ (see fig.\ \myref{fig:uang-level-weak} near $\xxo$).
\end{proposition}
\begin{proof}
  1. Let $r=\vlen{\xx}$ be distance from $\xxo$. 
  By assumption in Theorem \myref{th:weak-inf-univalent} $\stf$ is $\Cinf$ in a neighbourhood of $\xxo$ in $\cDomxy$ where it satisfies the equation
  \begin{alignat}{9} 0 = L\stf \defeq -\Aco:\nabx^2\stf \topref{eq:interior-stf}{=} -\Lap_\xx\stf + O(r^2):\nabx^2\stf  \myeqlabel{eq:Lap-O}  \end{alignat}
  as well as $\stf=0$ on $\sett{\ye=0}$ and on $\sett{\xe=\hw(\ye)}$. 
  Using that $\hw$ was also assumed $\Cinf$ we change to coordinates
  \begin{alignat}{9} \xi \defeq \hw(\ye)-\xe \geq 0 \csep  \eta \defeq \ye \geq 0 \quad, \end{alignat}
  with $\rho=\vlen{(\xi,\eta)}$ comparable to $r$.
  $\hw_{\ye}=O(\rho)$ near $\ye=0$, so
  \begin{alignat}{9} \stf_{\xe\xe}=\stf_{\xi\xi} \csep
  \stf_{\ye\ye}=\stf_{\eta\eta}+\hw_{\ye\ye}(\eta)\stf_\xi+O(\rho):\nabla_{(\xi,\eta)}^2\stf \quad; \end{alignat}
  the equation \myeqref{eq:Lap-O} turns into
  \begin{alignat}{9} 0 = -\Lap_{(\xi,\eta)}\stf - \hw_{\ye\ye}(\eta)\stf_\xi + O(\rho):\nabla_{(\xi,\eta)}^2\stf 
  . \end{alignat}
  For any $\eps>0$ the ansatz
  \begin{alignat}{9} \substf = 2\rho^\eps\xi\eta = \rho^{2+\eps}\sin(2\vartheta) \csep \vartheta=\arctan(\eta/\xi) \end{alignat}
  yields by $\substf_\xi=O(\rho^{1+\eps})$ and $\hw_{yy}=O(1)$ and $\nabla^2\substf=O(\rho^\eps)$ the result
  \begin{alignat}{9} L\substf = -\big( \subeq{ (2+\eps)^2-2^2 }{=(4+\eps)\eps} \big) \rho^\eps \sin(2\vartheta) + O(\rho^{1+\eps}) , \end{alignat}
  so in the sector $\sett{\xi,\eta\geq 0,~0<\rho\leq \delta'}$ for $\delta'>0$ sufficiently small we have a subsolution,
  meaning $L\substf\leq 0$ in the interior and
  obviously $\substf=0$ on the $\xi=0$ and $\eta=0$ boundaries including $(\xi,\eta)=0$. 

  The third boundary of the sector is the $\rho=\delta'$ quartercircle. 
  By \myeqref{eq:normal-deriv-stf-walls} the inner normal derivatives of $\stf$ on axis and wedge are positive,
  and they are the inward tangential derivatives of the quartercircle at the endpoints.
  Combined with $\substf=\stf=0$ in the endpoints and $\substf$ smooth on the closed quartercircle
  we find that $\stf\geq \iota \substf$ near the endpoints for $\iota >0$ sufficiently small,
  and on the entire compact quartercircle if $\iota >0$ is chosen sufficiently smaller.
  Since $\stf$ is a solution and $\iota \substf$ a subsolution,
  the comparison principle for $L$ shows $\stf\geq \iota \substf$ in the entire sector, in particular near the origin.
  Specializing to $\vartheta=\pi /4$ where $\substf=\rho^{2+\eps}$,
  we may 
  return to $\xe,\ye$ coordinates, finding a sequence of points $(\xx_j)\conv\xxo$ and a constant $\gamma>0$ so that
  \begin{alignat}{9} \stf(\xx_j) \geq  \gamma \vlen{\xx_j}^{2+\eps} \quad\text{for all $j$.} \myeqlabel{eq:stf-lowerbound-stag} \end{alignat}
  
  2. We assumed in Theorem \myref{th:weak-inf-univalent} that $\stf$ is $\Cinf$ at $\xxo$.
  Taking tangential derivatives of the slip condition $\stf=0$ along axis and wedge we obtain $\stf_\xe=0=\stf_\eta$ in $\xxo$. 
  However, if $\nabla^2\stf$ was also zero, then $\stf=O(r^3)$ would contradict \myeqref{eq:stf-lowerbound-stag}.
  Hence $\nabla^2\stf$ is nonzero.
  The equation \myeqref{eq:Lap-O} requires that the degree 2 part of the Taylor series of $\stf$ in $\xxo$ is harmonic and satisfies $\stf=0$ on axis and wedge,
  so necessarily it is a multiple of $\xe\ye$. $\stf$ was shown to be positive in $\Domxy$ (Proposition \myref{prop:stf-positive-level-sets}),
  in particular in the $\xe<0<\ye$ quadrant near $\xxo$, so we must have
  \begin{alignat}{9} \stf = - a_2\xe\ye + O(r^3) \end{alignat}
  for some constant $a_2>0$, with all $O$ terms still $\Cinf$ of course, and analytic away from the wedge (Remark \myref{rem:analyticity}).
  But then
  \begin{alignat}{9} \vv = -\dens^{-1}\nabla^\perp\stf
  &=
  \dens^{-1}a_2(-\xe,\ye) + O(r^2)
  \\\impl\quad
  \vang = \arctan\frac{\vy}{\vx} &= \arctan\frac{\ye}{-\xe} + O(r),   \quad\text{as well as}
  \\
  \nabx\vang &= \frac{(\ye,-\xe)}{r^2} + O(1) \neq 0  
  \end{alignat}
  showing the level sets are non-singular, with the one for $\vang_+\in\boi{0}{\pi/2}$ approaching near $\xxo$ the tangent with angle $\pi-\vang_+$ to the positive horizontal axis as $\xxo$ is approached (see $\xxo$ in figure \myref{fig:uang-level-weak}). 
\end{proof}

\begin{proposition}
  \mylabel{prop:angle-levels-wedge}%
  $\vang$ is strictly decreasing along $\Wedgexy$ from $\pi/2$ at $\ye\dnconv0$ to $\turnai$ at $\ye\conv\infty$.
  For $\vang_+\in\boi{\turnai}{\pi/2}$ the level set $L_+=\selset{\xx\in\Domxy}{\vang(\xx)=\vang_+}$ approaches $\Wedgexy$ in exactly one point,
  near which it is a simple analytic curve with non-tangential $\Cinf$ extension to the point.
  If $\vang_+\in\boi{0}{\turnai}$ then $L_+$ is separated from $\Wedgexy$. 
\end{proposition}\begin{proof}
  \newcommand{\vangW}{\vang^{\Wedge}}
  Let $(\cos\vangW,\sin\vangW)$ be the upward unit tangent to $\Wedgexy$,
  with continuous $\vangW$ fixed by $\vangW\conv\pi/2$ as $\ye\dnconv 0$, so $\vangW\conv\turnai$ as $\ye\upconv+\infty$ \myeqref{eq:hw-slope}.
  By definition before \myeqref{eq:vang-axis} $\vang$ is also continuous, 
  by the slip condition $\vang-\vangW$ is a constant multiple of $\pi$ along $\Wedgexy$,
  and by Proposition \myref{prop:angle-levels-stag} we have $\vang\conv\pi/2$ at $\Wedgexy$ as $\ye\dnconv 0$, so the constant is zero: $\vang=\vangW$.
  So the level set is separated from $\Wedgexy$ if $\vang_+\in\boi{0}{\turnai}$, 
  but for $\boi{\turnai}{\pi/2}$ approaches exactly one point of $\Wedgexy$. 
  By \emph{locally uniform} convexity \myeqref{eq:wedge-convex} the tangential part of $\nabx\vangW$ is nonzero,
  so by the implicit function theorem the extended level set is at that point a non-tangential simple $\Cinf$ curve that is analytic in its interior because $\vang$ is. 
\end{proof}

\begin{proposition}\mylabel{prop:velocity-rightward}%
  \begin{alignat}{9}  \vang &\in \cli{-\pi/2}{\pi/2}  \quad\text{in $\cDomxy\setdiff\sett{\xxo}$.} \myeqlabel{eq:vang-interval}\end{alignat} 
\end{proposition}\begin{proof}
  $\vang$ does not have extrema in the interior $\eDomxy$, by Proposition \myref{prop:vang-interior-extrema}.
  $\vang$ was fixed to $0$ on $\Axis\union\sett{\xxn}$ \myeqref{eq:vang-axis}.
  It has values in $\boi{\vang_w}{\pi/2}$ at $\Wedgexy$ by Proposition \myref{prop:angle-levels-wedge}.
  At $\cShockxy$ the values of $\vang$ are contained in $\cli{-\turnamax}{\turnamax}$ (the upper shock is connected and $\vang=0$ in its endpoint $\xxn$).
  Finally, its values on small circles around $\xxo$ from $\Axisxy$ to $\Wedgexy$ are limited to $\cli{0}{\pi/2+\eps}$ (Proposition \myref{prop:angle-levels-stag}),
  with $\eps>0$ shrinking to $0$ as the circle radius decreases.
  This proves the result in the upper halfplane, and the lower halfplane follow by the assumption of mirror symmetry \myeqref{eq:mirror-symmetry}. 
\end{proof}

\begin{proposition}
  \mylabel{prop:vang-extrema-shock}%
  If $\vang$ attains a local extremum at the shock, then the value is $\pm\turnamax$ and the shock curvature is negative in the extremum, in the sense $\hs_{\ye\ye}<0$.
\end{proposition}
(We emphasize that here and elsewhere ``local extremum at'' means
that it is an extremum over $V\isect\cDomxy$, \emph{not} just over $V\isect\Shockxy$, for $V$ a disk of sufficiently small radius $\delta>0$ centered in the boundary point.
\begin{proof}
  Using the notation \myeqref{eq:vang-arctan} we choose
  Cartesian coordinates $(\nu,\tau)$ so that \emph{in} the extremum point (but not necessarily elsewhere)
  $\partial_\nu$ is unit normal \emph{inward} to the downstream region. $\partial_\tau$ is 90 degree counterclockwise from $\partial_\nu$. 

  At the downstream side of the shock the values of $\vv$ are on the $\vv$-shock polar (section \myref{section:shockpolars}). 
  Consider a local extremum of $\vang$ 
  at the shock.
  Assume $\vang\neq\pm\turnamax$ there, then $\vang$ is locally a parameter for the $\vv$-shock polar, 
  so the first-order condition $\vang_\tau=0$ for an extremum implies $0=\vv_\tau$, i.e.\ $0=\vpot_{\tau\tau}=\vpot_{\nu\tau}$.
  But then the equation $0=\Aco:\nabx^2\vpot$ shows $\vpot_{\nu\nu}=0$ as well (using $\Aco^{\nu\nu}>0$ by positive definiteness), so $\vang_\nu=0$.
  But $0\neq\vang_\nu$ is required for a local extremum at a smooth boundary, by the Hopf lemma \cite{hopf-lemma}.

  Therefore $\vang$ can only have value $\pm\turnamax$ in the extremum (these are the extrema on the shock polar,
  i.e.\ largest and smallest possible \emph{at} the shock, but not elsewhere of course;
  $\vang\in\cli{-\pi/2}{\pi/2}$ by \myeqref{eq:vang-interval} does not permit adding multiples of $2\pi$). 

  Assume we have a local extremum with $\vang=+\turnamax$, necessarily a maximum. Then
  \begin{alignat}{9} \vang_\tau = 0 ,\end{alignat}
  and the Hopf lemma for the elliptic PDE \myeqref{eq:vang-pde} for $\vang$ requires that
  \begin{alignat}{9} \vang_\nu < 0 \myeqlabel{eq:vang-nu-negative} \end{alignat}
  (recall $\partial_\nu$ is \emph{inward}). Besides 
  \begin{alignat}{9} \vpot_\tau>0 \myeqlabel{eq:vpot-tau-positive} \end{alignat}
  (this is obvious from a graphical sketch, but here is a formal proof:
  $\vang=+\turnamax\in\boi0{\pi/2}$ means $\vy>0$;
  by \myeqref{eq:shockn-downstream} $\shockn=(\vvi-\vv)/\vlen{\vvi-\vv}$
  is a downstream-pointing unit shock normal, $\vvi=(\vxi,0)$, and by definition of $\partial_\tau$ we have $\vpot_\tau=\vv\dotp\shockn^\perp$ 
  which evaluates to $\vxi\vy/\vlen{\vvi-\vv}>0$). By definition of $\partial_\nu$ we have
  \begin{alignat}{9} \vpot_{i\nu}-\vpot_\nu = (\vvi-\vv)\dotp\shockn > 0 \myeqlabel{eq:vpot-nu} \end{alignat}
  for admissible shocks.

  Recalling \myeqref{eq:vang-3x3} we argue by Cramer's rule that
  \begin{alignat}{9} \vpot_{\tau\tau} &= \frac{\Det\begin{bmatrix}  \Acosym^{\nu\nu} & 2\Acosym^{\nu\tau}  & 0 \\  -\vpot_\tau/\vlen{\nabla\vpot}^2 & \vpot_\nu/\vlen{\nabla\vpot}^2  & \vang_\nu \\  0           & -\vpot_\tau/\vlen{\nabla\vpot}^2 & 0 \end{bmatrix}}{ \nabla\vpot\dotp\Aco\nabla\vpot/\vlen{\nabla\vpot}^4 }
  \notag\\&= \subeq{ \frac{ \Acosym^{\nu\nu} }{ \nabla\vpot\dotp\Aco \nabla\vpot / \vlen{\nabla\vpot}^2 } }{>0} \subeq{ \vpot_\tau }{ >0 } \subeq{ \vang_\nu }{<0} < 0 . \myeqlabel{eq:vpot-tau-tau}\end{alignat}
  using \myeqref{eq:vang-nu-negative}, \myeqref{eq:vpot-tau-positive} and $\Aco$ positive definite.

  Let $\partial_t$ be the tangential derivative \emph{everywhere} on the shock (as opposed to $\partial_\tau$ which need not be tangential except in the extremum point).
  Let $\curv$ be curvature in the extremum point, oriented to be positive when $\hs_{yy}>0$;
  we use that $\partial_\tau^2=\partial_t^2-\curv \partial_\nu$:
  $\vpot_i$ is linear and $(\tau,\nu)$ Cartesian, so $\vpot_{i\tau\tau}=0$, showing
  \begin{alignat}{9} 0 \topref{eq:vpot-tau-tau}{<} \subeq{\vpot_{i\tau\tau}}{=0} - \vpot_{\tau\tau} = \subeq{(\vpoti-\vpot)_{tt}}{=0} - \curv(\subeq{\vpot_{i\nu}-\vpot_\nu}{>0})
  \quad\toprefb{eq:vpot-nu}{eq:shock-jmp-vpot}{\impl}\quad \curv < 0 . \end{alignat} 
  
  For $\vang=-\turnamax$ we simply use mirror symmetry (or repeat the argument above with some inequalities reversed).
\end{proof}

\newcommand{\shocknca}{\shockn^{c*}}

\begin{proposition}
  \mylabel{prop:vang-nonnegative}%
  $\vang$ does not attain negative values in the above-axis domain $\cDomxy\setdiff\sett{\xxo}$.
\end{proposition}\begin{proof}
  Assume it does.
  By \myeqref{eq:vang-axis} $\vang$ is nonnegative at $\Axisxy$,
  by Proposition \myref{prop:angle-levels-wedge} also at $\Wedgexy$ and post-shock upper infinity,
  and by \myref{prop:angle-levels-stag} near $\xxo$.
  By Proposition \myref{prop:vang-interior-extrema} It does not have extrema in the interior $\Domxy$.
  So it must attain its minimum $\vang_-<0$ over $\cDomxy\setdiff\sett{\xxo}$ at the shock $\Shockxy$.

  By Proposition \myref{prop:vang-extrema-shock} the negative $\vang_-$ can only be $\turnamax$, and shock curvature is \emph{negative} in any $\vang$ extremum at the shock, in particular in any point where $\vang=-\turnamax$. However, at these points $\vv$ can only be the mirror image $\vvca$ of $\vvc$ on the lower half-polar, attained only by the shock normal $\shocknca=(\vvi-\vvca)/\vlen{\vvi-\vvca}$. The set of these points is compact and separated from $\xxn$ (where $\vv\conv\vvn\neq\vvca$) and upper infinity (where $\vv\conv\vvw\neq\vvca$),
  so there is at least one such point with minimal $\xx\dotp\shocknca$. But since $\shocknca$ points strictly right- and upward, the rest of the shock is on the \emph{right} side of the line through that point with normal $\shocknca$, so the shock curvature is \emph{nonnegative} there --- contradiction.
\end{proof}

\newcommand{\Lc}{L_c}%
\newcommand{\Lca}{L_{c*}}%

\begin{proof}[Proof of Theorem \myref{th:weak-inf-univalent}]\ \par
  I. Consider the ``critical'' level set $\Lc=\selset{\xx\in\cDomxy\setdiff\sett{\xxo}}{\vang=\turnamax}$. 
  By Proposition \myref{prop:levels-tree} $\Lc\isect\Domxy$ is a disjoint union of locally finite trees.
  
  (a) Assume some tree in $\Lc\isect\Domxy$ has a vertex. By Proposition \myref{prop:levels-atpoint} the degree is $\geq 4$,
  so there are at least four simple paths, with pairwise disjoint interiors, to $\bdry\Domxy$ or upper post-shock infinity.

  $\Lc$ is separated from $\cAxisxy\setdiff\sett{\xxo}$, where $\vang=0<\turnamax$ \myeqref{eq:vang-axis},
  and from infinity where $\vang=\turnai<\turnamax$. 
    By Proposition \myref{prop:angle-levels-wedge} $\vang$ strictly decreases along $\Wedgexy$,
  from $\pi/2$ in $\xxo$ along the segment $\Wedgexy^2$
  to $\turnamax$ in a point $\xxcW$
  and then along the segment $\Wedgexy^w$
  to $\turnai$ at upper infinity (fig.\ \myref{fig:uang-level-weak}); exactly one $\Lc$ edge approaches $\Wedgexy$, ending in $\xxcW$.
  By Proposition \myref{prop:angle-levels-stag} exactly one $\Lc$ edge approaches $\xxo$.
  Hence at least two simple paths approach $\Shockxy$.
  By Proposition \myref{prop:shock-interior-analytic} $\Shockxy$ is an analytic curve, with $\stf$ and $\vang$ analytic in $\Shockxy\union\Domxy$, 
  and by Remark \myref{rem:single-endpoint-boundary} each edge that approaches the shock meets it in exactly one point
  (otherwise $\vang=\turnamax$ everywhere on the shock, with obvious contradictions).
  The two simple paths meet the shock in two points, perhaps identical, and form with the shock segment between them a closed simple loop.
  By the Jordan curve theorem the loop encloses a nonempty interior with compact closure $K$. 
  By continuity (we excluded the edge to $\xxo$) the nonconstant $\vang$ attains a minimum $<\turnamax$ over $K$,
  by Proposition \myref{prop:vang-interior-extrema} necessarily on $\bdry K$, hence at $\Shockxy$,
  but by Proposition \myref{prop:vang-extrema-shock} that can only be $-\turnamax$, contradicting Proposition \myref{prop:vang-nonnegative}.

  Hence there are no interior vertices; $\Lc$ consists of edges with pairwise disjoint interiors, with endpoints $\xxcW$ or $\xxo$ or on $\Shockxy$.
  Each edge has at most one endpoint on $\Shockxy$, by the same $\vang$ minimum argument,
  and because $\xxcW$ and $\xxo$ each belong to exactly one edge (Propositions \myref{prop:angle-levels-wedge} and \myref{prop:angle-levels-stag}), 
  there cannot be more than two edges. 

  (b) (i) If $\Lc$ has only one edge, then its endpoints are $\xxo,\xxcW$; any $\Lc$ point on $\Shockxy$ must be isolated.

  (ii) If $\Lc$ has two edges, $\Lc^s$ starting in $\xxo$ and $\Lc^w$ starting in $\xxcW$,
  then their other endpoints must be on $\Shockxy$.
  If the endpoints are distinct, then the $\Shockxy$ segment between them has points with $\vang<\turnamax$
  (otherwise $\vang=\turnamax$ along the entire segment (contradiction), as $\turnamax$ is by definition the largest possible $\vang$ at the shock).
  But the segment $\Wedgexy^2$ from $\xxo$ to $\xxcW$ has $\vang>\turnamax$,
  so edges and segments form another simple loop $\bdry K$ with $\vang$ minimum $<\turnamax$ at $\Shockxy$, contradiction.
  Hence $\Lc^s,\Lc^w$ have the same $\Shockxy$ endpoint $\xxcS$, which is the unique \emph{non-isolated} $\Lc$ point of $\Shockxy$. 
  
  (c) \emph{Isolated} critical points on $\cShockxy$ are local maxima of $\vang$,
  so by Proposition \myref{prop:vang-extrema-shock} the shock curvature is negative there:
  on the nearby \emph{upper} side slope $d\ye/d\xe$ is higher, closer to a normal shock, so the shock is \emph{strong}-type there; 
  on the nearby \emph{lower} side slope is lower and the shock is \emph{weak}-type. \\
  (i) The shock is normal on the symmetry axis, in particular strong-type, so the lowest of all critical points has \emph{strong}-type on its lower side, hence cannot be isolated.\\
  (ii) The shock is weak-type near infinity, so the highest of all critical points has \emph{weak}-type on its upper side, hence cannot be isolated either. \\
  But by (b) there is at most one non-isolated critical point.
  So $\xxcS$ is the \emph{only} critical point on $\Shockxy$. 
  The shock segment $\Shockxy^s$ from $\xxn$ to $\xxcS$ is strong-type,
  the segment $\Shockxy^w$ from $\xxcS$ to infinity is weak-type.

  II. $\cDomxy$ is split by $\Lc$ into three open regions. \\
  (a) The first region $\Domxy^s$ is enclosed by $\Axisxy$, $\Lc^s$ and $\Shockxy^s$.
  On $\Shockxy^s$ we have $\vang\in\boi{0}{\turnamax}$,
  by Proposition \myref{prop:angle-levels-stag} we have $\vang\in\boi{0}{\turnamax}$ for any level edge approaching $\xxo$ through $\Domxy^s$ (below $\Lc^s$ above $\Axisxy$),
  and $\vang=0$ on $\cAxisxy\setdiff\sett{\xxo}$,
  so by absence of interior extrema (Proposition \myref{prop:vang-interior-extrema}) necessarily $\vang\in\boi{0}{\turnamax}$ in $\Domxy^s$.
  
  In $\Domxy^s$ the $\vang$-level set $L_+$ for value $\vang_+\in\boi{0}{\turnamax}$ is a single edge from $\xxo$ to $\Shockxy^2$, by arguments similar to the $\Lc$ case:
  isolated $L_+$ points on $\Shockxy^s$ are local extrema, contradicting Propositions \myref{prop:vang-extrema-shock} and \myref{prop:vang-nonnegative},
  and if any of the $L_+$ trees had a vertex, then it would have no simple paths to $\cAxisxy$, at most one to $\xxo$, at least three to $\Shockxy^s$, 
  any two forming with some segment of $\Shockxy^s$ a simple loop $\bdry K$ yielding a minimum $<\vang_+$ at $\Shockxy^s$, contradiction.
  The same argument rules out single edges with both endpoints on $\Shockxy^s$. 
  So each $L_+$ edge has exactly one endpoint in $\xxo$, but by Proposition \myref{prop:angle-levels-stag} there is exactly one such edge.
  Hence each $L_+$ is one single edge from $\xxo$ to $\Shockxy^s$. 

  (b) The second region $\Domxy^2$ is enclosed by $\Lc^s$, $\Lc^w$ and $\Wedgexy^2$.
  Along the latter $\vang\in\boi{\turnamax}{\pi/2}$ is strictly monotone, so local extrema are impossible,
  and now-familiar arguments show $\vang\in\boi{\turnamax}{\pi/2}$ in $\Domxy^2$,
  with each level set a single edge from $\xxo$ to $\Wedgexy^2$. 

  (c) The third open region $\Domxy^w$ is unbounded, bordered by $\Shockxy^w$, $\Lc^w$, $\Wedgexy^w$ and upper post-shock infinity.
  $\vang\in\boi{\turnai}{\turnamax}$ at $\Wedgexy^w$ and on the weak-type $\Shockxy^w$, while $\vang\conv\turnai$ at infinity, 
  so $\vang\in\boi{\turnai}{\turnamax}$ in $\Domxy^w$, with each level set a single edge from $\Wedgexy^w$ to $\Shockxy^w$. 

  III. (a)
  $\vlen{\vv}$ is strictly monotone along every single $\vang$-level edge, because otherwise in some interior point the tangential derivatives
  $\vang_t$ and $\vlen{\vv}_t$ are both zero, i.e.\ $\vec t\dotp\nabx^2\vpot\,\vv^\perp=0=\vec t\dotp\nabx^2\vpot\,\vv$,
  showing by Lemma \myref{lemma:two-zero-three}(i) that $\nabx^2\vpot=0$, so $\nabx\vang=0$,
  but that is a vertex (Proposition \myref{prop:levels-atpoint}) --- contradiction.

  (b) It remains to determine the \emph{direction} of strict monotoniticity:
  along the $\vv$-shock polar $\vlen{\vv}$ is strictly increasing from normal to vanishing (section \myref{section:shockpolars});
  $\vlen{\vv}>\vlen{\vvc}$ on the weak-type $\Shockxy^w$, but $\vlen{\vv}<\vlen{\vvc}$ on strong-type $\Shockxy^s$.
  So in $\Domxy^s$ along each $\vang$-level edge $\vlen{\vv}$ is strictly increasing from $0$ in $\xxo$ to $<\vlen{\vvc}$ on $\Shockxy^s$.
  The levels $L_+$ for $\vang_+\upconv\turnamax$ converge to the boundary $\Lc^s$
  where continuity implies \emph{non-strict} increase of $\vlen{\vv}$ from $\xxo$ to $\xxcS$,
  but we already knew strict \emph{monotonicity}, combined we have strictly increasing.
  Passing to the other side $\Domxy^2$, again continuity shows strict increase along every level edge from lower to upper endpoint on $\Wedgexy^2$,
  again by limit on the other boundary $\Lc^s$ from $\xxcS$ to $\Wedgexy^s$,
  and then on the other side $\Domxy^w$ along every level edge from $>\vlen{\vvc}$ at $\Shockxy^w$ to a larger value on $\Wedgexy^w$.

  (c) We have shown $\vlen{\vv}$ is strictly monotone along $\Axisxy\union\sett{\xxo}$ where $\vang=0$,
  along $\Lc$ where $\vang=\turnamax$,
  and along each $\vang$-level set in $\Domxy^2\union\Wedgexy^2$ where $\vang>\turnamax$,
  in $\Domxy^s\union\Shockxy^s$ where $0<\vang<\turnamax$ and $\vlen{\vv}<\vlen{\vvc}$,
  and in $\Domxy^w\union\Shockxy^w$ where $0<\turnai<\vang<\turnamax$ but $\vlen{\vv}>\vlen{\vvc}$.
  For $\ye<0$ we use mirror symmetry with $\vang\leftarrow -\vang$, and we have $\vv=0$ only in $\xxo$.
  So since $\xx\mapsto\vv$ is injective on each of these subsets and the value sets are pairwise disjoint, it is injective on the entire domain $\ecDomxy$.
\end{proof}

\begin{figure}
  \centerline{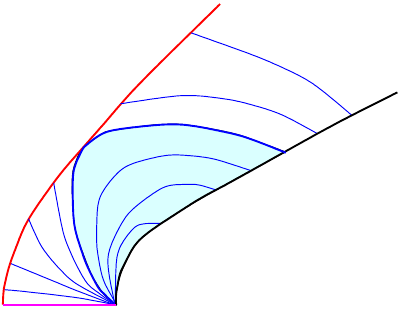}
  \caption{Velocity angle $\measuredangle\vv$ level sets for weak-at-infinity solutions}
  \myfiglabel{fig:uang-level-weak}
\end{figure}

\begin{remark}
  \mylabel{rem:strong-at-inf}
  For a \emph{strong}-at-infinity solution, step I(c)(ii) of the preceding proof is no longer valid, but can be adapted to show the following alternative:\\
  1. either there is no critical point on the shock at all, with critical level set from $\xxo$ to the wedge-critical point $\xxcW$ through the interior,\\
  2. or there is a lower non-isolated critical point as in the weak-at-infinity case, with one edge to $\xxo$ and another to $\xxcW$,
  but there is an \emph{upper} \emph{isolated} critical point. \\
  In the first case the flow is univalent, so existence proofs for the former case are straightforward by modifications of the following half of this paper. 
  But in the second case the flow is \emph{not} univalent because $\xx\mapsto\vv$ attains the value $\vvc$ twice on the shock. 
  Proving that such flows do exist requires entirely different methods, so both cases will be discussed in future work.
\end{remark}

\mysection{Equations in the $\pq$ plane}\mylabel{section:pq-equations}

In the second half of the paper we prove Theorem \myref{th:weak-inf-existence}, existence for weak-at-infinity flows onto a blunt convex wedge.
The first half proves that any such flow is univalent, with injective $\cDomxy\ni\xx\mapsto\jj$.
So we can construct solutions in the $\jj$ plane, then map back to $\xx$, both of which are nontrivial steps.
The benefit of working in the $\jj$ plane is that the shock is a \emph{fixed} boundary, a segment of the mass flux polar.
It is also convenient that the interior equation becomes linear, since the literature about corners for such equations is a bit more detailed. 

Since $\xx\mapsto\vv$ is also injective, the reader may wonder why the standard approach of working in the $\vv$ plane is avoided.
This is because the slip condition for $\vpot$ takes on rather awkward forms in these coordinates, a technical obstacle that can be bypassed by working with $\stf$ rather than $\vpot$.
Using $\stf$ comes with other technical problems, but as we show all of them can be solved.

We change to new coordinates
\begin{alignat}{9} \pq = (\pqvar^1,\pqvar^2) = \nabx\stf = \jj^\perp . \myeqlabel{eq:pq-def}\end{alignat}
For deriving the $\pq$ plane equations assume that $\xx\mapsto\pq$ is a diffeomorphism (say $\Cinf$), which of course requires proof later on.
We define a new potential
\begin{alignat}{9} \pqpot(\pq)  \defeq   \pq^T\xx(\pq) - \stf(\xx(\pq)). \myeqlabel{eq:pqpot-def}\end{alignat}
Its derivative 
\begin{alignat}{9} \nabp\pqpot = \xx + \subeq{\nabp\xx^T\pq - \nabp\xx^T \subeq{\nabx\stf}{=\pq}}{=0}  \myeqlabel{eq:xx-nabp-pqpot}\end{alignat} 
gives the position $\xx$ where the mass flux $\jj=-\pq^\perp$ is located.
\begin{alignat}{9} \nabp^2\pqpot = \xx_\pq = \pq_\xx^{-1} = (\nabx^2\stf)^{-1} = (\Det\nabx^2\stf)^{-1} \adjt\nabx^2\stf , \end{alignat}
so the interior PDE
\begin{alignat}{9} 0 = \Aco:\nabx^2\stf \end{alignat}
is equivalent to
\begin{alignat}{9} 0 = \Aco:\adjt\nabp^2\pqpot &= \Acosym^{11}\pd 2^2\pqpot - \Acosym^{12}\pd 1\pd 2\pqpot + \Acosym^{22}\pd 1^2\pqpot
\\&= \adjt\Aco:\nabp^2\pqpot. \myeqlabel{eq:interior-pq}\end{alignat}
Given a solution $\pqpot$, the original $\stf$ is recovered by \myeqref{eq:pqpot-def} in the form
\begin{alignat}{9}  \stf(\xx) \defeq  \pq(\xx)^T\xx - \pqpot(\pq(\xx)), \myeqlabel{eq:stf-pqpot}\end{alignat}
where $\pq$ is found by solving \myeqref{eq:xx-nabp-pqpot} $\nabp\pqpot(\pq)=\xx$.

At slip condition walls $\stf$ is constant \myeqref{eq:slip}, but in our case we fixed the constant to $0$:
\begin{alignat}{9} 0 \topref{eq:slip-stf}{=} \stf \topref{eq:stf-pqpot}{=} \pq\dotp\nabp\pqpot - \pqpot .  \myeqlabel{eq:slipcond-pq} \end{alignat}

$\pq$ being located on $\Shockpq$ is one of the two shock conditions.
The other one is \myeqref{eq:shock-jmp-stf} $\stf=\stfi$ expressed via \myeqref{eq:stf-pqpot} as
\begin{alignat}{9} \stfi = \pq\dotp\nabp\pqpot - \pqpot .\end{alignat}
$\stfi$ is linear with gradient $\pqi={\jji}^\perp=(0,\densi\vxi)$, and we fixed (\myeqref{eq:shock-jmp-stf} and above) $\stfi=0=\stf$ at the normal point $\xxn$, so
$\stfi = \pqiy\ye = \pqi\dotp\xx$, giving the final form
\begin{alignat}{9} 0 = (\pq-\pqi)\dotp\nabp \pqpot - \pqpot  . \myeqlabel{eq:polarcond} \end{alignat}

Conversely, the conditions \myeqref{eq:polarcond} and $\pq\in\Shockpq$ imply the original two shock conditions \myeqref{eq:shock-jmp-stf} and \myeqref{eq:shock-ut}:\\
(i) \myeqref{eq:polarcond} with definition \myeqref{eq:stf-pqpot} shows \myeqref{eq:shock-jmp-stf} $\stf=\stfi$,
which implies by tangential derivative that $\nabla\stf-\nabla\stfi = (\jj-\jji)^\perp $ is a shock normal;\\
(ii) $\Shockpq$ is part of the $\pq$-shock polar, which is the rotated $\jj$-polar, defined by \myeqref{eq:jj-polar-eqn} as the set of points $\jj$ where 
\begin{alignat}{9} 0 = (\dens^{-1}\jj-\densi^{-1}\jji)\dotp(\jj-\jji),  \end{alignat}
so using (i) this implies
\begin{alignat}{9} 0 = (\vv-\vvi)\dotp\shockn^\perp \end{alignat} 
which is \myeqref{eq:shock-ut}.

\begin{figure}%
  \centerline{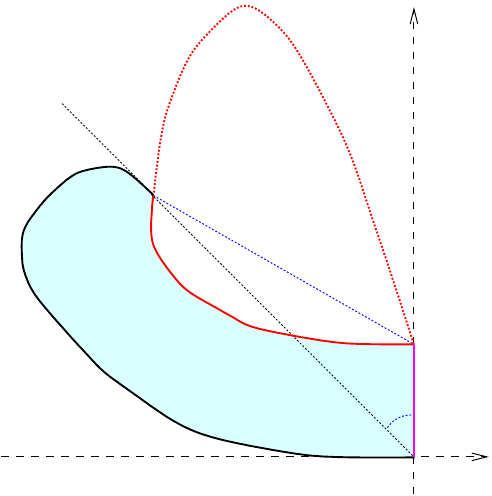}
  \caption{$\pq$ plane domain; $\Wedgepq$ can be chosen quite arbitrarily subject to some conditions.}
  \myfiglabel{fig:pqdomainsymm}
\end{figure}

As before $\Shockpq,\Wedgepq,\Axispq$ refer to the interior of curves, excluding endpoints; see figure \myref{fig:pqdomainsymm}:
\begin{definition}\mylabel{def:domain}
  Let $\cShockpq$ be the segment of the subsonic $\pq$ shock polar from some point $\pqw$ on its weak-type segment to the normal point $\pqin$.
  $\cAxispq$ is the line segment from $\pqo$ to $\pqin$.
  $\cWedgepq$ is permitted to be any curve \\
  (i) through $\sett{\vlen\pq<\jmax,\pqvar^2<0<\pqvar^1}$ \\
  (ii) that is $\Cinf$ non-degenerate including the endpoints, \\
  (iii) analytic in the interior, \\
  (iv) converges to $\pqo$ and to $\pqw$ at the endpoints, \\
  (v) has inward unit tangents $(0,-1)$ at $\pqo$ and $\pqw/\vlen\pqw$ at $\pqw$, \\
  (vi) does not intersect itself or $\Shockpq\union\sett\pqi\union\Axispq$, \\
  (vii) with angle $\measuredangle\pq$ locally uniformly monotone on the interior. 

  By the Jordan curve theorem the simple closed curve $\cShockpq\union\cWedgepq\union\cAxispq$ encloses a single bounded connected component $\Dompq$. By compactness
  \begin{alignat}{9} \sup_{\cDompq}\Mach<1. \myeqlabel{eq:MachboundDom}\end{alignat}
\end{definition}

\begin{remark}
  For a flow as in Theorem \myref{th:weak-inf-univalent}, l.u.\ convex wedge corresponds to the angle condition (vii).
  Condition (vii) is equivalent to the vector $\pq$ not being tangential to $\Wedgepq$ in any interior point $\pq$,
  which means the slip condition \myeqref{eq:slipcond-pq} $\pq\dotp\nabp\pqpot=\pqpot$ is nowhere degenerate oblique. 
  In two dimensions this is less important than in three or higher, but it permits citing freely the large majority of literature that is stated for arbitrary dimension,
  for example \cite[Theorem 6.26]{gilbarg-trudinger}.
  If the condition is relaxed, univalent flows for some non-convex wedges can be obtained by essentially the same proof. 
\end{remark}

\mysection{Maximum principles for $\dex\dotp\nabp\pqpot$}

Let $\dex$ be a constant nonzero vector.
Since the interior PDE \myeqref{eq:interior-pq}
\begin{alignat}{9} 0 = \adjt\Aco:\nabp^2\pqpot = \Acosym^{22}\pd 1^2\pqpot - 2\Acosym^{12}\pd 1 \pd 2\pqpot + \subeq{\Acosym^{11}}{>0}\pd 2^2\pqpot \myeqlabel{eq:interior-pq-repeated} \end{alignat} 
is symmetric under rotations around the origin, we may assume $\dex=(1,0)$. 
$\pd 1$ gives
\begin{alignat}{9} 0 = \adjt\Aco:\nabp^2\pd 1\pqpot + \pd 1(\adjt\Aco):\nabp^2\pqpot. \myeqlabel{eq:interior-pd1} \end{alignat}
In the second term  the $\pd 2^2\pqpot$ component of $\nabp^2\pqpot$ can be eliminated using \myeqref{eq:interior-pq-repeated} 
since $\Acosym^{11}>0$ by positive definiteness, so \myeqref{eq:interior-pd1} reduces to
\begin{alignat}{9} 0 = \adjt\Aco:\nabp^2\pd 1\pqpot + \vec{\tilde{a}}\dotp\nabp\pd 1\pqpot \end{alignat}
for some function $\vec{\tilde a}$ of $\pq$ (depending on $\pqpot$). Returning to original coordinates we have 
\begin{alignat}{9} 0 &= \adjt\Aco:\nabp^2(\dex\dotp\nabp)\pqpot
+ \vec a\dotp\nabp(\dex\dotp\nabp)\pqpot \myeqlabel{eq:dex-nabp-pqpot-pde}\end{alignat}
for some rotated function $\vec a$ of $\pq$.

\begin{proposition}\mylabel{prop:dex-pqpot-extrema-interior}%
  For $\dex$ a constant nonzero vector,
  $\dex\dotp\nabp\pqpot$ does not have extrema in the interior $\Dompq$ unless it is constant. 
\end{proposition}\begin{proof}
  The strong maximum principle applied to the elliptic PDE \myeqref{eq:dex-nabp-pqpot-pde} shows $\dex\dotp\nabp\pqpot$ cannot have extrema in the interior unless it is constant.
\end{proof}

\begin{proposition}\mylabel{prop:dex-pqpot-extrema-wall}%
  If $\dex$ is a constant nonzero vector, then non-constant $\dex\dotp\nabp\pqpot$ cannot attain extrema at $\Wedgepq$ unless it is a multiple of $\dex$. 
\end{proposition}\begin{proof}
  In an extremum the tangential derivative is zero:
  \begin{alignat}{9} 0 = \ttp\dotp\nabp^2\pqpot\dex .\end{alignat}
  The boundary condition \myeqref{eq:slipcond-pq} has tangential derivative
  \begin{alignat}{9} 0 = \ttp\dotp\nabp(\pq\dotp\nabp\pqpot-\pqpot) = \ttp\dotp\nabp^2\pqpot\pq + \subeq{\ttp\dotp\nabp\pqpot - \ttp\dotp\nabp\pqpot}{=0} . \end{alignat}
  Comparing the last two equations, Lemma \myref{lemma:two-zero-three} (i) shows that if $\pq$ is not a multiple of $\dex$, then $\nabp^2\pqpot=0$ in the extremum.
  But then the normal derivative $\nnp\dotp\nabp^2\pqpot\dex$ is also zero, so the Hopf lemma applied to the elliptic PDE \myeqref{eq:dex-nabp-pqpot-pde} shows $\dex\dotp\nabp\pqpot$ is constant.
\end{proof}

\begin{proposition}\mylabel{prop:dex-pqpot-extrema-shock}%
  If $\dex$ is a constant nonzero vector, then nonconstant $\dex\dotp\nabp\pqpot$ cannot attain extrema at $\Shockpq$ except where $\pq-\pqi$ is a multiple of $\dex$. 
\end{proposition}\begin{proof}
  The boundary condition $\pq\dotp\nabp$ in \myeqref{eq:slipcond-pq} is replaced by $(\pq-\pqi)\dotp\nabp$ in \myeqref{eq:polarcond}, so in the proof of
  Proposition \myref{prop:dex-pqpot-extrema-wall} replace $\pq$ by $\pq-\pqi$. 
\end{proof}

\begin{proposition}\mylabel{prop:dex-vertical-ok}%
  (i) $\pq$ is not a multiple of $\dex=(0,1)$ on $\Wedgepq$.\\
  (ii) $\pq-\pqi$ is not a multiple of $\dex=(0,1)$ on $\Shockpq$. \\
  (iii) $\ye=\pd2\pqpot$ cannot have local extrema except in $\pqw$ or at $\cAxispq$.
\end{proposition}\begin{proof}
  (i) 
  $\Wedgepq$ was chosen in Definition \myref{def:domain}(v)+(vii) so that $\measuredangle\pq$ is strictly increasing from $\pi/2$ (vertical) at $\pq\conv\pqo$ to $\measuredangle\pqw\in\boi{\pi/2}{\pi}$ in $\pqw$,
  so $\pq$ is not vertical in the interior.
  (ii) $\pqi$ is vertical, so by section \myref{section:shockpolars} $\pq-\pqi$ is not vertical on the interior of $\Shockpq$.
  ($\pqi\parallel\pq-\pqi$ means $\vvi\parallel\vv-\vvi$ so $\vvi\parallel\vv$ which happens only for normal shocks.)
  (iii) $\dex$ satisfies by (i) the conditions of Proposition \myref{prop:dex-pqpot-extrema-wall} and by (ii) those of Proposition \myref{prop:dex-pqpot-extrema-shock},
  combined with Proposition \myref{prop:dex-pqpot-extrema-interior} ruling out extrema of $\pd2\pqpot$ anywhere except in $\pqw$ or on $\cAxispq$. 
\end{proof}

\mysection{Shift solution}

\begin{proposition}
  \mylabel{prop:pqpotshift}%
  For $\dex=(1,0)$, 
  \begin{alignat}{9}   \pqpotshift(\pq)  \defeq  \dex\dotp(\pq-\pqi)   \myeqlabel{eq:pqpotshift}   \end{alignat}
  is a solution of \myeqref{eq:interior-pq} and \myeqref{eq:polarcond} and \myeqref{eq:slipcond-pq}.
\end{proposition}\begin{proof}
  The interior pde \myeqref{eq:interior-pq} $\adjt\Aco:\nabd^2\pqpotshift=0$ is clearly satisfied by any linear $\pqpot$.

  Polar condition \myeqref{eq:polarcond}:
  \begin{alignat}{9} (\pq-\pqi)\dotp\nabp\pqpotshift-\pqpotshift = (\pq-\pqi)\dotp\dex-\dex\dotp(\pq-\pqi) = 0 \end{alignat}
  for any $\dex$. 

  Slip condition \myeqref{eq:slipcond-pq}:
  \begin{alignat}{9} 0 = \pq\dotp\nabp\pqpotshift-\pqpotshift = \pq\dotp\dex - \dex\dotp(\pq-\pqi) = \dex\dotp\pqi \end{alignat}
  which is true if and only if $\pqi=\jji^\perp=(0,\densi\vxi)$ is perpendicular to $\dex$, i.e.\ if $\dex$ is a multiple of $(1,0)$. 
\end{proof}
$\pqpotshift$ is not a useful solution of its own, since $\pq\mapsto\nabp\pqpotshift=\dex$ is constant rather than injective, but adding multiples of $\pqpotshift$ to other solutions $\pqpot$ performs horizontal rigid translation in the $\xx$ plane.

\mysection{Dilation}\mylabel{section:dilation}

\newcommand{\Rdil}{z}

To compute the asymptotics near each corner it is easier to change to \defm{dilated} coordinates $\dpq$ where the principal part of the operator is a Laplacian.
The derivative matrix $\dpq_{\pq}$ of the coordinate transformation performs a dilation in direction of $\vv$ and a different dilation in direction of $\vv^\perp$, as shown below.

In new coordinates $\dpq$ the principal part has coefficient matrix a multiple of 
\begin{alignat}{9} \dpq_\pq \Aco^{-1} \dpq_\pq^T \end{alignat}
For this to be a multiple of the identity matrix it is sufficient that $\dpq_\pq$ is a multiple of a square root of the symmetric matrix $\Aco=\idm-(\vv/\csnd)^2$,
which has eigenvectors $\vv$ (or $\jj$ or $\pq^\perp$) for eigenvalue $1-\Mach^2$ and $\vv^\perp$ (or $\jj^\perp$ or $\pq$) for eigenvalue $1$.
A suitable choice is 
\begin{alignat}{9} \dpq &\defeq e^\Rdil \pq \myeqlabel{eq:dilation}, 
  \csep \Rdil &\defeq \int_0^{\vlenq{\pq}} \frac{ (1-\Mach^2)^{-1/2}-1 }{\vlenq{\pq}} d(\vlenq{\pq}). \end{alignat}
$(1-\Mach^2)^{\pm 1/2}$ is analytic, in particular bounded, on the compact set $\cDompq$ where $\sup_{\cDompq}\Mach<1$ \myeqref{eq:MachboundDom};
the integrand and $\Rdil$ are analytic there as well, even at $\pq=0$ where the numerator in \myeqref{eq:dilation} is also zero, with expansion
\begin{alignat}{9} \dens\csnd &= \dens_0\csnd_0 + O(\vlenq\pq) \end{alignat}
from the discussion after \myeqref{eq:dHpm-ddens-atjj}, so
\begin{alignat}{9} \Mach^2 &= \frac{\vlenq{\pq}}{(\dens\csnd)^2} = \frac{\vlenq\pq}{(\densstag\csndstag)^2} + O(\vlenq{\pq}) , 
\\ \frac{ (1-\Mach^2)^{-1/2}-1 }{\vlenq{\pq}} &= \frac{ \half\Mach^2 +O(\Mach^4) }{\vlenq{\pq}} = \frac1{2(\densstag\csndstag)^2} + O(\vlenq{\pq}),
\\ \Rdil &= 0 + \frac1{2(\densstag\csndstag)^2} \vlenq\pq + O(\vlen{\pq}^4) .\end{alignat}
The derivative of $\pq$ is
\begin{alignat}{9}
\dpq_\pq &= \frac1{2(\densstag\csndstag)^2} \idm
\quad \text{in $\pq=0$,}
\\
\dpq_{\pq} &= e^\Rdil (\frac{d\Rdil}{d\vlenq{\pq}}\pq^2+\idm)
\\&= e^\Rdil \Big((\vlenq{\pq}\frac{d\Rdil}{d\vlenq{\pq}}+1)(\frac{\pq}{\vlen\pq})^2+(\frac{\pq^\perp}{\vlen\pq})^2\Big) 
\\&=  e^\Rdil \Big( (1-\Mach^2)^{-1/2} (\frac{\pq}{\vlen\pq})^2 + (\frac{\pq^\perp}{\vlen\pq})^2 \Big)
\quad \text{in $\pq\neq 0$,}
\myeqlabel{eq:Ddil-explicit}
\end{alignat}
with eigenvectors $\pq/\vlen\pq,\pq^\perp/\vlen\pq$ for eigenvalues $1,(1-\Mach^2)^{1/2}$ times the same scalar $e^{\Rdil}(1-\Mach^2)^{-1/2}$, 
so $\dpq_{\pq}$ is indeed a square root of $\Aco$. $\dpq_\pq$ and its inverse are analytic and therefore uniformly bounded on $\cDompq$.

The equation \myeqref{eq:interior-pq} in $\Dompq$, which is homogeneous so that nonzero scalar factors can be removed, transforms to 
\begin{alignat}{9} 0 = (\subeq{\Lap_{\dpq} + \firstcoeff\dotp\nabd}{=\inop(1)}) \pqpot \quad\text{in $\dDompq$,} \myeqlabel{eq:inop1} \end{alignat}
here $\dDompq$ is the image of $\Dompq$ under $\pq\mapsto\dpq$, while $\firstcoeff$ is some analytic function that does not matter for our purposes.
The boundary conditions will be discussed in detail later.

\mysection{Corner exponents}\mylabel{section:corner-exponents}

To discuss existence of solutions, as well as asymptotics in the corners, we use methods for linear elliptic PDE and systems from \cite{nazarov-plamenevskii}, which follows the 
approach of Kondrat'ev \cite{kondratev} and the literature following it \cite{mazya-plamenevskii,mazya-plamenevskii-asy}.
This approach takes a Mellin transform in the radial variable, which is a Fourier-Laplace transform in the log-radius,
producing a family of operators parametrized by the dual variable $\lambda$.
The principal parts of perators and domain are constant-coefficient elliptic operators in cones. 
The transformed family of principal parts, called ``operator pencil'' in \cite{nazarov-plamenevskii},
has ``eigenfunctions'' corresponding to homogeneous solutions for ``eigenvalues'' whose imaginary part $\Im\lambda$ corresponds to rate of decay/growth in the corner. 
We refer to this part as \defm{corner exponent}. 
In our case the principal part of the interior operator is the simplest possible one, a planar Laplacian, but at the boundary we have
oblique derivative conditions with variable coefficients and singularities, so detailed calculations are needed.

For $\cori=o,i,w$ let $\rcori$ be distance to to the corner $\dpqcori$ in $\dpq$ coordinates.
$\rcori$ is uniformly positive and bounded outside any neighbourhood of $\dpqcori$, by continuity and compactness of $\cdDompq$. 
$\rcori$ is also comparable to the distance from $\pq$ to $\pq^\cori$ because $\dpq_\pq$ and its inverse $\pq_\dpq$ are bounded (section \myref{section:dilation}).

The equations to solve are (compare section \myref{section:pq-equations})
\begin{alignat}{9}  0  &=  ( \pq\dotp\nabp - 1 ) \pqpot  &&\qquad\text{on $\Axispq$,}   \myeqlabel{eq:bcond-Axis}
\\  0  &=  ( \pq\dotp\nabp - 1 ) \pqpot  &&\qquad\text{on $\Wedgepq$,}  \myeqlabel{eq:bcond-Wedge}
\\  0  &=  ( (\pq-\pqi)\dotp\nabp - 1 ) \pqpot   &&\qquad\text{on $\Shockpq$,}   \myeqlabel{eq:bcond-Shock}
\\  0  &=  ( \Lapd + \firstcoeff\dotp\nabd ) \pqpot  &&\qquad\text{in $\dDompq$,} \myeqlabel{eq:bcond-Dom}\end{alignat}
where using \myeqref{eq:inop1} the interior operator is already in dilated coordinates. 
We may equivalently multiply each equation by some scalar function nonzero in $\cdDompq\setdiff\sett{\dpqo,\dpqi,\dpqw}$:
\begin{alignat}{9}  0  &=  \ri ( \pq\dotp\nabp - 1 ) \pqpot   &&\quad\text{on $\Axispq$,}  \myeqlabel{eq:bcond-Axis-scaled}
\\  0  &=  \rw ( \pq\dotp\nabp - 1 ) \pqpot   &&\quad\text{on $\Wedgepq$,}  \myeqlabel{eq:bcond-Wedge-scaled}
\\  0  &=  \rw ( (\pq-\pqi)\dotp\nabp - 1 ) \pqpot   &&\quad\text{on $\Shockpq$,}  \myeqlabel{eq:bcond-Shock-scaled}
\\  0  &=  (\ro\ri\rw)^2 ( \Lapd+\firstcoeff\dotp\nabd ) \pqpot   &&\quad\text{in $\Dompq$.} \myeqlabel{eq:inopr}\end{alignat}
Now our operators match the setting of \cite[p.\ 99]{nazarov-plamenevskii}, where the radial derivatives come scaled as $\radpr$. 

To calculate the Fredholm index of our problem we use the standard technique of constructing an entire path of Fredholm operators, parametrized by $\lpr\in\cli01$ and continuous in the operator norm topology, with $\lpr=1$ our desired operator, while $\lpr=0$ is a simpler operator corresponding to the Laplace-Neumann problem.
This choice, and the following choice of connecting path, is quite arbitrary and a matter of convenience.
We extend the interior operator \myeqref{eq:inopr} to
\begin{alignat}{9}  \inop(\lpr)  &\defeq  (\rw\ro\ri)^2 ( \Lapd + \lpr \firstcoeff\dotp\nabd ) . \myeqlabel{eq:inoplpr} \end{alignat}
We write the boundary operators as
\begin{alignat}{9} \bdop(\lpr) \defeq \dbb(\lpr) \dotp \nabd + \bzero(\lpr)   \myeqlabel{eq:bdop-lpr}  \end{alignat}
($\dbb(\lpr),\bzero(\lpr)$ functions of $\dpq$). For $\lpr$ from $0$ to $\half$ we change the derivative direction $\dbb$ from normal to the desired value, while $\lpr$ from $\half$ to $1$ is for changing the order $0$ coefficient $\bzero$ from $0$ to its desired value:
\begin{alignat}{9} \bzero(\lpr) &\defeq \begin{cases}
    0 , & 0\leq \lpr\leq \half , \\
    (2\lpr-1)\bzero(1) , & \half \leq \lpr\leq 1,
    \myeqlabel{eq:bzero-lpr}
\end{cases} \end{alignat}
with $\bzero(1)$ determined by our desired problem;
\begin{alignat}{9}  \dbb(\lpr)  \defeq  \vlen{\dbb(1)} \, \big( \cos \dbcang(\lpr) \dts + \sin \dbcang(\lpr) \dts^\perp \big)  \myeqlabel{eq:dbb-lpr}  \end{alignat}
where $\dts$ is the \emph{counterclockwise} tangent to the smooth parts of $\bdry\dDompq$ so that $\dts^\perp$ is the \emph{inward} unit normal
(fig.\ \myref{fig:cornernotationalt}), 
$\dbcang$ the \emph{counterclockwise} angle from $\dts$ to $\dbb$, with
furem{
  keep $\vlen{\dbb_B(\lpr)}$ b/c $O(r)$ sometimes, cannot interp to $O(1)$, instead need to interp $0$th order part from $1$ to $0$.}
\begin{alignat}{9}  \dbcang(0)  \defeq  \pi /2,  \myeqlabel{eq:kappa-lprzero}\end{alignat}
which corresponds to a Neumann condition, and 
\begin{alignat}{9} \dbcang(\lpr) \defeq \begin{cases}
  (1-2\lpr)\dbcang(0)+2\lpr\dbcang(1) ,& 0\leq \lpr\leq \half , \\
  \dbcang(1)                   ,& \half \leq \lpr\leq 1. \end{cases} \myeqlabel{eq:kappa-lpr}\end{alignat}
Again $\dbcang(1)$ is determined by our desired problem (up to a multiple of $2\pi $, which will be chosen below on each smooth boundary segment). We combine interior and boundary operators into
\begin{alignat}{9}  \allop(\lpr)  =  (\inop(\lpr),\bdop(\lpr)) .\myeqlabel{eq:allop-def}\end{alignat}

\begin{figure}
  \centerline{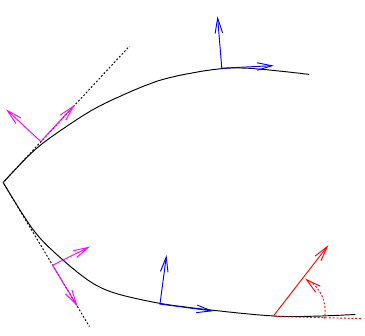}
  \caption{Notation in some corner of $\dDompq$}
  \myfiglabel{fig:cornernotationalt}
\end{figure}

Consider polar coordinates $(\rad,\omegav)$ in each corner. The domain converges to the sector $\omegav\in\boi{\omegalo}{\omegahi}$ at the corner (fig.\ \myref{fig:cornernotationalt}).

The first-order part of the interior operator \myeqref{eq:inoplpr} is $\rad^2\lpr\firstcoeff\dotp\nabd$ for $\lpr\firstcoeff$ bounded; this does not belong to the principal part of the operator because $\rad^2$, as opposed to $\rad$, decays too fast. The principal part is entirely from $\Lapd$, in polar coordinates 
\begin{alignat}{9} 0 = ((\radpr)^2 + \pom^2)\pqpot  \quad\text{for $\rad>0,\omegav\in\boi{\omegalo}{\omegahi}$.} \end{alignat}
The solutions are linear combinations of
\begin{alignat}{9}& \rad^{\critex}\sin(\critex\omegav),\rad^{\critex}\cos(\critex\omegav) \quad\text{for $\critex\in\RR\setdiff\sett0$,}  \myeqlabel{eq:comb-q-notzero}
\\& 1,\ln\rad,\omegav \quad\text{corresponding to $\critex=0$.}   \myeqlabel{eq:comb-q-zero} \end{alignat}
$\critex$ is what \cite{nazarov-plamenevskii} call ``imaginary part'' of the ``points of the spectrum of the pencil'' of ``principal homogeneous parts of the operators'';
for obvious reasions we prefer the term \defm{corner exponent}. After boundary conditions are applied, the space spanned by \myeqref{eq:comb-q-notzero} or \myeqref{eq:comb-q-zero} reduces to a finite dimension called \defm{multiplicity} of $\critex$ here and in \cite{nazarov-plamenevskii}. 

$\dts$ was chosen counterclockwise around the domain, pointing away from the corner on the $\omegav\conv\omegazero$ boundary but towards on the other side
(fig.\ \myref{fig:cornernotationalt}),
so
\begin{alignat}{9} \rad(\dts\dotp\nabd,\dts^\perp\dotp\nabd) \conv \begin{cases}
  (\radpr,\pom) , & \text{on the $\omegav\conv\omegazero$ side,} \\
  (-\radpr,-\pom) , & \text{on the $\omegav\conv\omegaone$ side.}
\end{cases} \myeqlabel{eq:dtsnabconv} \end{alignat}

Consider $0\leq \lpr\leq \half $, so $\bzero(\lpr)=0$, so the boundary conditions \myeqref{eq:bdop-lpr}--\myeqref{eq:kappa-lpr} have principal parts equivalent to 
\begin{alignat}{9} 0 = \vlen{\dbb(1)} (\cos \dbcangjside,\sin \dbcangjside)\dotp( \radpr  , \pom  )\pqpot  \quad\text{on $\omegav=\omegajside$} \qquad(\jside=0,1), \myeqlabel{eq:bcondomj}\end{alignat}
where $\dbcangjside$ is the limit of $\dbcang$ on the $\omegav\conv\omegajside$ side of the corner.
\begin{enumerate}
\item
  For \defm{corner exponent} $\critex\neq 0$ write the combinations of \myeqref{eq:comb-q-notzero} as $\pqpot=\rad^\critex\sin(\critex\omegav-\gamma)$ for constant $\gamma$, then
  \begin{alignat}{9}  0  =  \critex \rad^\critex \big(  \subeq{  \cos \dbcang_j \sin(\critex\omegajside-\gamma) + \sin \dbcangjside \cos(\critex\omegajside-\gamma)  }{  =\sin(\critex\omegajside-\gamma+\dbcangjside)  }  \big) \quad(\jside=0,1), \end{alignat}
  which is equivalent to
  \begin{alignat}{9} \critex\omegajside-\gamma+\dbcangjside = \pi \ice_\jside \quad(\jside=0,1) \end{alignat}
  for integer $\ice_\jside$. With $\ice=\ice_1-\ice_0$ 
  we get the solutions
  \begin{alignat}{9} \critex &= \frac{ \pi \ice - (\dbcanga1-\dbcanga0) }{\omegaone-\omegazero} \quad\text{for any $\ice \in \ZZ$,} \myeqlabel{eq:critex-Dkappa}
  \\ \gamma &= \dbcanga0 + \critex\omegazero , \end{alignat}
  where adding multiples of $\pi $ to $\gamma$ would only change the sign of 
  \begin{alignat}{9} \pqpot = r^\critex \sin(\critex(\omegav-\omegazero)-\dbcanga0). \myeqlabel{eq:critex-pqpot}\end{alignat}
\item
  For $\critex=0$: among \myeqref{eq:comb-q-zero} $\pqpot=1$ is always a solution of \myeqref{eq:bcondomj}; the ansatz $\pqpot=c_t\ln\rad+c_\omegav\omegav$ for constants $c_t,c_\omegav$ leads to
  \begin{alignat}{9} 0 = c_t \cos \dbcangjside + c_\omegav \sin \dbcangjside \quad(\jside=0,1),\end{alignat}
  which has nontrivial solutions $(c_t,c_\omegav)$ if and only if $(\cos \dbcanga0,\sin \dbcanga0)$ and $(\cos \dbcanga1,\sin \dbcanga1)$ are parallel. 
  If so, then 
  \begin{alignat}{9} \pqpot = (-\sin \dbcanga0,\cos\dbcanga0)\dotp(\ln\rad,\omegav) \myeqlabel{eq:pqpot-log}\end{alignat}
  spans together with $\pqpot=1$ a two-dimensional space of homogeneous solutions. 
\end{enumerate}

\noindent The case $\half \leq \lpr\leq 1$ is discussed in context for each corner.

\subsection{Stagnation corner $\pqo$}

Consider the stagnation corner $\pq=\pqo$ (fig.\ \myref{fig:pqdomainsymm}).
At the corner the outward unit tangents converge to $(-1,0)$ along $\Wedgepq$, to $(0,1)$ along $\Axispq$. 
Those are eigenvectors of $\dpq_\pq$ in $\pq=\pqo$ \myeqref{eq:Ddil-explicit}, so the interior corner angle is $\pi /2$ before and after dilation;
we pick $\omegazero\defeq0$ for $\dAxispq$ so that $\omegav\conv\omegaone\defeq\pi/2$ along $\dWedgepq$. 

Near the $\veczero$ corner $\dpq=\dpq_{\pq}\pq+O(\rad^2)$,
so the boundary conditions \myeqref{eq:bcond-Axis-scaled} and \myeqref{eq:bcond-Wedge-scaled} are equivalent to
\begin{alignat}{9} 0
&= \pqpot_\pq \pq - \pqpot
\\&= \pqpot_\dpq \dpq_\pq \pq - \pqpot
\\&= \pqpot_\dpq ( \dpq + O(\rad^2)) - \pqpot
\\&= \pqpot_\dpq \big( \subeq{\vlen{\dpq}}{=\rad} \subeq{\frac{\dpq}{\vlen{\dpq}}}{(*)} + O(\rad^2)\big) - \pqpot
\\&= \radpr \pqpot + \pqpot_\dpq O(\rad^2) - \pqpot \end{alignat}
because $(*)$ converges to an outward unit tangent on either boundary.
Clearly the principal part is
\begin{alignat}{9} 0 &= \radpr \pqpot - \pqpot .\end{alignat}
In \myeqref{eq:bdop-lpr} and \myeqref{eq:dbb-lpr} the corresponding values are (recall \myeqref{eq:dtsnabconv})
\begin{alignat}{9} \dbcanga0(1)=0=\dbcanga1(1) \csep \bzero_0(1) = -1 \csep \bzero_1(1)=+1 , \end{alignat}
so for $\lpr\in\cli0{\half}$ we get on both sides
\begin{alignat}{9}  \dbcangjside(\lpr) \toprefb{eq:kappa-lpr}{eq:kappa-lprzero}{=} (\half -\lpr)\pi   \quad(\jside=0,1). \end{alignat}
So the discussion after \myeqref{eq:bcondomj} shows the homogeneous solutions are multiples of 
\begin{alignat}{9} \pqpot \topref{eq:critex-pqpot}{=} \rad^\critex\sin\big(\critex\omegav-(\half -\lpr)\pi \big) \csep  \critex \topref{eq:critex-Dkappa}{=} 2\ice, \myeqlabel{eq:stag-hom-nonzero-critex}\end{alignat}
for nonzero $\critex$; for zero $\critex$ the homogeneous solutions are linear combinations of
\begin{alignat}{9} \pqpot = 1 \csep \pqpot \topref{eq:pqpot-log}{=} \cos(\lpr\pi )\ln \rad+\sin(\lpr\pi )\omegav .\myeqlabel{eq:stag-hom-zero-critex}\end{alignat}
For $\lpr\in\cli{\half}{1}$ by \myeqref{eq:bzero-lpr} the principal parts of the boundary conditions on both sides are equivalent to
\begin{alignat}{9} 0 &= \big( \radpr  - (2\lpr-1) \big)\pqpot . \end{alignat}
The homogeneous solutions are the $\lpr=\half$ case of \myeqref{eq:stag-hom-nonzero-critex},
\begin{alignat}{9} \pqpot = \rad^{\critex} \sin(\critex\omegav) \csep \critex = 2\ice \csep \ice\in\ZZ\setdiff\sett0, \end{alignat} 
along with linear combinations of
\begin{alignat}{9} \pqpot = r^{2\lpr-1}\sin((2\lpr-1)\omegav) \csep \pqpot = r^{2\lpr-1}\cos((2\lpr-1)\omegav) \end{alignat}
for double exponent $\critex=2\lpr-1$. 
In summary, in the $\pqo$ corner the characteristic exponents (each counting $1$ towards multiplicity) are
\begin{alignat}{9}  \critex^o_{0\sigk}(\lpr)  &=  \begin{cases} 0 , & 0 \leq  \lpr \leq  \half  \\ 2\lpr-1 , & \half  \leq  \lpr \leq  1 \end{cases}  \quad\text{for $\sigk=\pm $,} 
\notag\\  \critex^o_{\ice\sigk}(\lpr)  &=  2\sigk\ice  \quad\text{for $\ice\in \NN$, $\sigk=\pm $.}  \myeqlabel{eq:critex-o} \end{alignat}

\subsection{Normal corner $\pqi$}

Similar to the stagnation corner $\pqo$,
as we approach the normal corner $\pqi=(0,\densi\vxi)$
the matrix $\dpq_\pq$ converges to a diagonal matrix \myeqref{eq:Ddil-explicit},
while the boundary tangents converge to eigenvectors $(-1,0)$ along $\Shockpq$ and $(0,1)$ along $\Axispq$, 
so the corner angle is $\pi /2$ before and after dilation.

Still similar to the $\pqo$ corner,
near the $\pqi$ corner $\dpq-\dpqi=\dpq_\pq(\pq-\pqi)+O(\rad^2)$,
so the boundary condition \myeqref{eq:bcond-Shock-scaled} on $\Shockpq$
is equivalent to
\begin{alignat}{9} 0 &= \pqpot_\pq(\pq-\pqi) - \pqpot
\\&= \pqpot_{\dpq}\dpq_\pq(\pq-\pqi) - \pqpot
\\&= \pqpot_{\dpq}\big(\dpq-\dpqi + O(\rad^2)\big) - \pqpot
\\&= \pqpot_{\dpq}\big(\subeq{\vlen{\dpq-\dpqi}}{=\rad} \subeq{\frac{\dpq-\dpqi}{\vlen{\dpq-\dpqi}}}{(*)} + O(\rad^2)\big) - \pqpot
\\&= \radpr\pqpot_{\dpq} + \pqpot_\dpq O(\rad^2) - \pqpot  \end{alignat}
because $(*)$ converges to an outward unit tangent. 
So the $\lpr=1$ principal part is
\begin{alignat}{9} (\radpr -1)\pqpot   = 0  \quad\text{on $\omegav=\omegazero\defeq 0$.} \myeqlabel{eq:bcond-normal-S} \end{alignat}
Approaching the normal corner $\pqi$ along the other boundary $\Axispq$, we argue $\pq=\pqi+O(\rad)$ is an exact tangent,
but it does not converge to zero any more as it did at the other endpoint,
so the boundary condition \myeqref{eq:bcond-Axis-scaled} is equivalent to
\begin{alignat}{9} 0 = (\vlen{\pqi}+O(r))\radpr \pqpot - O(r)\pqpot . \end{alignat}
Here the $0$th order term does not enter the $\lpr=1$ principal part, which is equivalent to
\begin{alignat}{9} \radpr \pqpot = 0  \quad\text{on $\omegav=\omegaone\defeq \pi /2$.} \quad\myeqlabel{eq:bcond-normal-X}\end{alignat}
For $0\leq \lpr\leq \half $ we have
\begin{alignat}{9} \dbcangjside(\lpr)=(\half -\lpr)\pi  \qquad(\jside=0,1) \end{alignat}
on both sides, so \myeqref{eq:critex-Dkappa} gives the same result as \myeqref{eq:stag-hom-nonzero-critex} and \myeqref{eq:stag-hom-zero-critex} in the stagnation corner $\pqo$:
\begin{alignat}{9} \pqpot \topref{eq:critex-pqpot}{=} \rad^\critex\sin\big(\critex\omegav-(\half -\lpr)\pi \big) \csep  \critex \topref{eq:critex-Dkappa}{=} 2\ice , 
\\ \pqpot = 1 \csep \pqpot \topref{eq:pqpot-log}{=} \cos(\lpr\pi )\ln \rad+\sin(\lpr\pi )\omegav \csep \critex=0. \end{alignat}
For $\half <\lpr\leq 1$ the principal parts of the boundary conditions are
\begin{alignat}{9} 0 &= (\radpr  - (2\lpr-1))\pqpot  \quad\text{on $\omegav=\omegazero=0$,}
\\ 0 &= \radpr \pqpot \quad\text{on $\omegav=\omegaone=\pi/2$,} \end{alignat}
with homogeneous solutions
\begin{alignat}{9} \pqpot = r^{\critex}\sin(\critex\omegav) \csep \critex = 2\ice \csep \ice \in \ZZ\setdiff\sett0 \end{alignat}
as well as
\begin{alignat}{9} \pqpot &= r^{\critex}\cos(\critex\omegav-\lpr \pi) \csep \critex = 2\lpr-1, 
\\ \pqpot &= \omegav \csep \critex = 0 . \end{alignat}
In summary, the characteristic exponents, each counting $1$ towards multiplicity, are
\begin{alignat}{9}  \critex^i_{0+}(\lpr)  &=  \begin{cases} 0 , & 0\leq \lpr\leq \half  , \\ 2\lpr-1 , & \half  \leq  \lpr \leq  1 , \end{cases}  
\notag\\  \critex^i_{0-}(\lpr)  &=  0  ,
\notag\\  \critex^i_{\ice \sigk}(\lpr)  &=  2\sigk\ice   \quad\text{for $\ice \in  \NN$, $\sigk=\pm $.}  \myeqlabel{eq:critex-n} \end{alignat}

\subsection{Weak corner $\pqw$}

\begin{figure}
  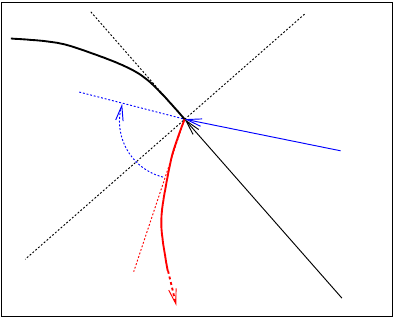\hfil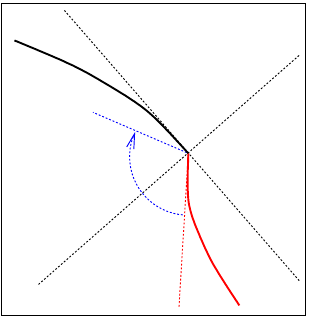
  \caption{Left: weak corner before dilation, right: after.
    Infinitesimally each of the four directions $\pm\pqw,\pm(\pqw)^\perp$ is mapped to itself, and each quadrant of directions enclosed is also mapped to itself.}
  \myfiglabel{fig:pqweakcorner}
\end{figure}

$\Shockpq$ was chosen (in Definition \myref{def:domain}) to be the polar segment from $\pqw$ to normal,
with endpoint $\pqw$ chosen to be on the \emph{weak-type subsonic} segment of the polar,
along which mass flux $\vlen{\pq}=\vlen{\jj}$ and angle $\measuredangle\pq=\vang+\pi/2$ are monotone (as recalled in section \myref{section:shockpolars}). 
More precisely: with $\partial_t$ the tangential derivative in direction from vanishing to normal,
\begin{alignat}{9}  0  &>  \partial_t\vlen{\pq}  =  \frac{ \pq\dotp \partial_t\pq }{ \vlen\pq }
\quad\text{and}\quad
0  <  \partial_t\measuredangle\pq  =  \frac{ \pq^\perp\dotp \partial_t\pq }{ \vlenq\pq }  \myeqlabel{eq:dtpq-ineqs}. \end{alignat}
In the corner $\pqw$ the outward unit tangent to $\Shockpq$ is precisely $\partial_t\pq / \vlen{\partial_t\pq}$,
so \myeqref{eq:dtpq-ineqs} shows it is somewhere in the interior of the quadrant counterclockwise from $\pqw^\perp$ to $-\pqw$.
Since $\pqw,\pqw^\perp$ are eigenvectors of $\dpq_\pq$ at $\pq=\pqw$ for positive eigenvalues $\lambda,\lambda^\perp$ \myeqref{eq:Ddil-explicit}, after dilation the outward unit tangent to $\dShockpq$ is still in the same quadrant,
at an angle
\begin{alignat}{9} \omegaow &\in \boi{\pi/2}{\pi}  \myeqlabel{eq:daw-interval} \end{alignat}
counterclockwise from the $\pqw$ direction. 

The outward tangent to the other boundary $\Wedgepq$ was chosen in Definition \myref{def:domain}(v) to be the eigenvector $\pqw$, 
so after dilation $\dWedgepq$ has the same outward tangent, which we assign to $\omegav=\omegazero\defeq 0$. 

The boundary condition \myeqref{eq:bcond-Wedge-scaled} is equivalent to
\begin{alignat}{9} 0
&=
\rw( \pqpot_\pq\pq - \pqpot )
\\&=
\rad\pqpot_\dpq\dpq_\pq\pq + O(\rad)\pqpot
\\&=
\rad\pqpot_\dpq\dpq_\pq (\pqw + O(\rad))  + O(\rad)\pqpot
\\&=
\rad\pqpot_\dpq (\lambda\pqw + O(\rad))  + O(\rad)\pqpot
\\&=
\subeq{\lambda\vlen{\pqw}}{\eqs 1}\, \rad \subeq{ \pqpot_\dpq \frac{\pqw}{\vlen{\pqw}} }{=\pr\pqpot+O(\rad)}  + \pqpot_\dpq O(\rad^2) + O(\rad)\pqpot,
\end{alignat} 
so the principal part is equivalent to
\begin{alignat}{9} 0 = \radpr\pqpot \quad\text{on $\omegav=\omegazero$,} \myeqlabel{eq:bcond-princ-Ww} \end{alignat}
corresponding to
\begin{alignat}{9} \dbcanga0(1)=0 \csep \bzero_0(1)=0 .\end{alignat}

Near $\pqw$ corner along $\Shockpq$ (see fig.\ \myref{fig:pqweakcorner}):
the $\pqw-\pqi$ direction is strictly inward to $\Dompq$ both on $\Shockpq$ and on $\Wedgepq$,
and this infinitesimal property is preserved by dilation, so $\dpq_\pq(\pqw-\pqi)$ is also strictly inward, at an angle
\begin{alignat}{9} -\dbcangow \in \boi{0}{\omegaone} \myeqlabel{eq:kappaSw} \end{alignat}
\emph{clockwise} from the $\omegav=\omegaone$ ray corresponding to $\dShockpq$. Boundary condition:
\begin{alignat}{9} 0 &\topref{eq:bcond-Shock-scaled}{=} \rw \big( \pqpot_\pq(\pq-\pqi) - \pqpot \big)
\\&= \rad\pqpot_\dpq\dpq_\pq\big(\pqw-\pqi+O(\rad)\big) + O(\rad)\pqpot 
\\&=
\subeq{\vlen{\dpq_\pq\dotp(\pqw-\pqi)}}{\approx\vlen{\dpqw-\dpqi}\eqs 1}
\,\rad\pqpot_\dpq\frac{\dpq_\pq(\pqw-\pqi)}{\vlen{\dpq_\pq\dotp(\pqw-\pqi)}} + \pqpot_\dpq O(\rad^2) + O(\rad)\pqpot,
\end{alignat}
so by definition of $\dbcangow$ the principal part is equivalent to 
\begin{alignat}{9} 0 = \begin{bmatrix}
  \rad\pqpot_{\rad}  \\
  \pqpot_{\omegav}
\end{bmatrix} \dotp \begin{bmatrix}
  \cos(-\dbcangow) \\
  \sin(-\dbcangow)
\end{bmatrix}
\quad\text{on $\omegav=\omegaone$.} \myeqlabel{eq:bcond-princ-Sw}\end{alignat}
corresponding to
\begin{alignat}{9} \dbcanga1(1) = \dbcangow  \csep \bzero_1(1)=0.  \end{alignat}

For $\lpr\in \cli0\half$ we have
\begin{alignat}{9}  \dbcanga0(\lpr) = (\half -\lpr)\pi    \csep  \dbcanga1(\lpr)=2\lpr\dbcangow+(\half -\lpr)\pi   , \end{alignat}
so \myeqref{eq:critex-Dkappa} provides corner exponents and homogeneous solutions
\begin{alignat}{9} \pqpot = \rad^\critex \sin\big(\critex\omegav-(\half -\lpr)\pi \big)  \csep  \critex = \frac{\ice \pi -2\lpr\dbcangow}{\omegaow} \csep \ice\in\ZZ.
\myeqlabel{eq:pqpot-weak-special}\end{alignat}
We estimate
\begin{alignat}{9} \abs{2\lpr\dbcangow} \overset{0\leq 2\lpr\leq 1}{\leq } \abs{\dbcangow}
\topref{eq:kappaSw}{<} \omegaone \topref{eq:daw-interval}{<} \pi \leq \abs{\ice\pi} \quad\text{when $\ice\neq 0$,} \end{alignat}
so $\critex\neq 0$ unless $\lpr=\ice=0$.
In that case (see \myeqref{eq:pqpot-log}) $\pqpot=\ln r$ is a particular solution with $\critex=0$. $\pqpot=1$ for $\critex=0$ is a particular solution in all cases.

For $\lpr$ from $\half$ to $1$ the boundary conditions do not change at all since there are no $0$th order parts. 

We summarize:
\begin{alignat}{9}  \critex^w_{0-}(\lpr)  &=  0, 
\myeqlabel{eq:critex-w-zerominus}
\\  \critex^w_{0+}(\lpr)  &=  \min\{2\lpr,1\}\cdot\frac{-\dbcangow}{\omegaow},
\\  \critex^w_{\ice \sigk}(\lpr)  &=  \critex^w_{0+} + \sigk\ice \frac{\pi }{\omegaow}  \qquad(\ice\in \NN,\sigk=\pm ). \end{alignat}
The exponents for $\ice\geq 1$ are less important because 
\begin{alignat}{9} \critex^w_{\ice-}(\lpr) &= \frac{-\min\{2\lpr,1\}\dbcangow-\ice \pi }{\omegaow} \topref{eq:kappaSw}{<} \frac{\omegaow-\ice \pi }{\omegaow} \topref{eq:daw-interval}{<} 0 ,
\\ \critex^w_{\ice+}(\lpr) &= \frac{-\min\{2\lpr,1\}\dbcangow+\ice \pi }{\omegaow} \topref{eq:kappaSw}{>} \frac{\ice \pi }{\omegaow} \topref{eq:daw-interval}{>} 1; \end{alignat}
combined with \myeqref{eq:critex-o} and \myeqref{eq:critex-n} we find there is a constant $\eps_0>0$ so that 
\begin{alignat}{9}  \critex^\cori_{\ice \sigk}(\lpr) &\begin{cases} < -\eps_0 , & \sigk\ice \leq  -1 \\ \in  \cli01 , &  \ice = 0 \\ > 1+\eps_0 , & \sigk\ice \geq  1 \end{cases}
\quad\text{for $\lpr\in \cli01,\cori\in\sett{i,o,w}$}. \myeqlabel{eq:eps0}\end{alignat}
We will choose corner weights in a small neighbourhood of $\cli01$, so \myeqref{eq:eps0} means only corner exponents for $\ice=0$ will affect Fredholmness and Fredholm indices. 
We give symbols to the most important homogeneous solutions at $\lpr=1$:
\begin{alignat}{9}  \pqpot^w_{0-}  &=  1 ,  \myeqlabel{eq:pqpot-w-zerominus}
\\  \pqpot^w_{0+}  &=  r^{\critex}\sin(\critex \omegav) \csep \critex=\critex^w_{0+}(1) = \frac{-\dbcangow}{\omegaow}
\toprefb{eq:kappaSw}{eq:daw-interval}{\in} \boi01.  \myeqlabel{eq:pqpot-w-zeroplus} \end{alignat}

\subsection{Leading weak-corner term asymptotics}

The homogeneous solution \myeqref{eq:pqpot-w-zeroplus} will be particularly important for us,
because it describes the leading correction to the straight shock near infinity. Let
\begin{alignat}{9} \partial_\tau = \frac{\dpqw}{\vlen\dpqw}\dotp\nabd \csep \partial_\nu = \frac{\dpqw^\perp}{\vlen\dpqw}\dotp\nabd \quad; \end{alignat}
$\partial_\tau$ is the derivative in $\dpqw$ direction, outward tangential to $\dWedgepq$ in $\dpqw$, and $\partial_\nu$ is 90 degree counterclockwise from $\partial_\tau$, 
so the relation to polar coordinates is
\begin{alignat}{9} \partial_\tau = \cos \omegav \, \pr  - \rad^{-1} \sin \omegav \, \pom  ,
\\ \partial_\nu = \sin \omegav \, \pr  + \rad^{-1} \cos \omegav \, \pom  . \myeqlabel{eq:dtau-dnu-def}\end{alignat}
Applied to \myeqref{eq:pqpot-w-zeroplus} we find, with $\critex=\critex^w_{0+}$ until \myeqref{eq:yv-particular-blowup},
\begin{alignat}{9} \begin{bmatrix}
  \partial_\tau  \\
  \partial_\nu
\end{bmatrix}\pqpot^w_{0+}
&= \critex\,\rad^{\critex-1}  \begin{bmatrix}
  \cos\omegav \sin(\critex\omegav) - \sin\omegav \cos(\critex\omegav)  \\
  \sin\omegav \sin(\critex\omegav) + \cos\omegav \cos(\critex\omegav)  \end{bmatrix}
\\&= \critex\,\rad^{\critex-1}  \begin{bmatrix}
  \sin((\critex-1)\omegav)  \\
  \cos((\critex-1)\omegav)  \end{bmatrix}.
\myeqlabel{eq:da-db-asy} \end{alignat}
$(\critex-1)\omegav$ is linear, with value zero in $\omegav=0$, but at $\omegav=\omegaow$
\begin{alignat}{9}
(\critex-1)\omegav
\topref{eq:pqpot-w-zeroplus}{=}
(\frac{-\dbcangow}{\omegaow}-1)\omegaow &= -\dbcangow-\omegaow
\\&\topref{eq:kappaSw}{\in} \boi{-\omegaow}{0} \topref{eq:daw-interval}{\subset } \boi{-\frac{\pi }{2}}{0}. \end{alignat}
and therefore for any $\omegav\in\cli{0}{\omegaow}$. Hence
\begin{alignat}{9} \cos((\critex-1)\omegav) &> 0    \quad\text{and}
\notag\\ \sin((\critex-1)\omega) &\leq 0     \qquad\text{in $\omegav\in\cli{0}{\omegaow}$.}
\myeqlabel{eq:cos-sin-nearw} \end{alignat}
Dilation $\pq\mapsto\dpq$ preserves direction \myeqref{eq:dilation}, so 
\begin{alignat}{9} \frac{\pqw}{\vlen\pqw} = \frac{\dpqw}{\vlen\dpqw}  \csep  \frac{\pqw^\perp}{\vlen\pqw} = \frac{\dpqw^\perp}{\vlen\dpqw}. \end{alignat}
These are also eigenvectors of $\dpq_\pq$ in $\pq=\pqw$ for positive eigenvalues $\lambda,\lambda^\perp $, so for $\pq$ \emph{near} $\pqw$
\begin{alignat}{9}  \dpq_\pq \frac{\pqw}{\vlen\pqw}  =  \lambda \frac{\dpqw}{\vlen\dpqw}  +  O(\rad)
\csep
\dpq_\pq \frac{\pqw^\perp}{\vlen\pqw}  =  \lambda^\perp  \frac{\dpqw^\perp}{\vlen\dpqw}  +  O(\rad).  \myeqlabel{eq:dir-Or} \end{alignat}
Since the direction of $\pqw$ is by choice at the counterclockwise angle $\turnai\in\boi{0}{\pi/2}$ from the upward vertical direction of $\pd 2$, we can argue 
\begin{alignat}{9} - \pd 2\pqpot^w_{0+} 
&=  \big( \sin\turnai\frac{\pqw^\perp}{\vlen\pqw} - \cos\turnai \frac{\pqw}{\vlen\pqw} \big) \dotp \nabp\pqpot^w_{0+} 
\\&\topref{eq:dir-Or}{=}  \big( \lambda^\perp \sin\turnai\frac{\dpqw^\perp}{\vlen\dpqw} - \lambda\cos\turnai \frac{\dpqw}{\vlen\dpqw} + O(\rad) \big) \dotp \nabd\pqpot^w_{0+} 
\\&\topref{eq:dtau-dnu-def}{=}  \big(\lambda^\perp \sin\turnai+O(\rad)\big) \partial_\nu \pqpot^w_{0+} - \big(\lambda\cos\turnai+O(\rad)\big) \partial_\tau\pqpot^w_{0+} \end{alignat}
\begin{alignat}{9}&\topref{eq:da-db-asy}{=}
\critex \rad^{\critex-1} \Big(
\supeq{\lambda^\perp \sin\turnai}{>0} \cos((\critex-1)\omegav)
- \subeq{\lambda\cos\turnai}{>0} \sin((\critex-1)\omegav)
+ O(\rad)
\Big).
\notag\end{alignat}
By \myeqref{eq:cos-sin-nearw} the $\sin,\cos$ terms are positive on $\omegav\in\cli{0}{\omegaow}$,
and by continuity and compactness they have positive minimum on $\omegav\in\cli{-\eps}{\omegaow+\eps}$ for some $\eps>0$.
For $\rad$ sufficiently small no other $\omegav$ occur, and the positive minimum dominates $O(\rad)$, so the big parenthesis is uniformly positive. 
By \myeqref{eq:pqpot-w-zeroplus} $\critex-1<0$, so $\critex\rad^{\critex-1}$ blows up and therefore
\begin{alignat}{9} -\pd 2\pqpot^w_{0+} &\conv  +\infty \quad\text{as $\pq\conv\pqw$ in $\cDompq$.}   \myeqlabel{eq:yv-particular-blowup} \end{alignat}

\mysection{Existence in the $\pq$ plane}\mylabel{section:existence-pq}

Having fixed the operators we define their domains and codomains based on the weighted $L^2$-Sobolev spaces defined in \cite[p.\ 63--65 and 97--99]{nazarov-plamenevskii}. With left-hand side our notation, right-hand side theirs: 
\begin{alignat}{9}
\pqo &= x^1, \\
\pqi &= x^2, \\
\pqw &= x^3, \\
\wtoi &= \sobl-n/2-\bt^1,\  n = 2 ,
\\&= \sobl-n/2-\bt^2     \myeqlabel{eq:wt-bt} \\
\wtw &= \sobl-n/2-\bt^3,\ T=3\\
\wtvec = (\wtoi,\wtw) , &\ \quad \btvec = (1,\bt^2,\bt^3), \end{alignat}
\begin{alignat}{9}
\dDompq &= \nzplDom ,\\
\domsob{\sobl}{\wtvec} &= \dvlb{\sobl}{\btvec}(\nzplDom), \\
\ransob{\sobl}{\wtvec} &= \rvlb{\sobl}{\btvec}(\nzplDom), \\
\allop &= \mathscr{A}. \end{alignat}
We also define cutoff functions $\cutcori\in\Cinf(\cDompq)$ that are $1$ near the corner $\pqcori$ but $0$ near any other corner.
Our \defm{corner weights} $\wt$, which are essential for Fredholm property and index, are better than $\bt$ for making transparent the decay and growth as $\rad=\vlen{\pq-\pqcori}\dnconv 0$:
functions $\cutcori\rad^\critex$ belong to $\domsob{\sobl}{\wtvec}$ if $\critex>\wtcori$, but not if $\critex<\wtcori$.

Our problem has analytic data and coefficients, so the solutions we expect are analytic, and the Sobolev exponent $\sobl$ has no effect on Fredholm index.
We only need it large enough for taking traces of enough derivatives to make the boundary conditions well-defined.
$\domsob\sobl\wtvec$ for fixed $\wtvec$ is decreasing as $\sobl$ grows \cite[(3.5.2), p.\ 63]{nazarov-plamenevskii}. 
We will use the convenient abbreviation
\begin{alignat}{9} \domsob\infty\wtvec \defeq \bigcap_{\sobl=1}^\infty \domsob\sobl\wtvec
\csep
\ransob\infty\wtvec \defeq \bigcap_{\sobl=1}^\infty \ransob\sobl\wtvec .\end{alignat}

\begin{proposition}
  \mylabel{prop:pqexi}%
  There exists a solution $\pqpot$,
  analytic in $\cDompq\setdiff\sett{\pqw,\pqi,\pqo})$ and $\Conce(\cDompq\setdiff\sett{\pqw})$,
  of the problem \myeqref{eq:bcond-Axis}--\myeqref{eq:bcond-Dom}
  such that
  \begin{alignat}{9}  \pd1\pqpot &= 0  \quad\text{in $\pq=\pqo$,}  \myeqlabel{eq:pqpot-nabp-zero-pqo} 
  \\  \pqpot = \pd2\pqpot &= 0  \quad\text{at $\cAxispq$,}  \myeqlabel{eq:pqpot-zero-Axis} \end{alignat}
  and such that for $\ye=\partial\pqpot/\partial\pqvar^2$ the limit $\yew\defeq\lim_{\pq\conv\pqw}\ye$ exists, and
  \begin{alignat}{9} \text{either}\quad  \yew &=  +\infty     \myeqlabel{eq:eitheryblowup}
  \\  \text{or}\quad  \pqpot \in \Conce(\cDompq) \csep \yew &\geq 0 .  \myeqlabel{eq:orynonnegative} \end{alignat}
\end{proposition}
\begin{proof}
  1. 
  $\allop(0)$ is the Laplace-Neumann operator.
  By \cite[section 6.1.4]{nazarov-plamenevskii} for sufficiently small $\eps>0$ we have Fredholm index
  \begin{alignat}{9} \find ( \allop(0) : \domsob\sobl{-\eps,-\eps} \rightarrow  \ransob\sobl{-\eps,-\eps} ) = +3  \myeqlabel{eq:find-dirneunegeps}  \end{alignat}
  because our domain has $Q=3$ corners.

  It is easy to check that $\cli01\ni\lpr\mapsto\allop(\lpr)$ is a continuous path in the $\domsob\sobl{-\eps,-\eps}\rightarrow \ransob\sobl{-\eps,-\eps}$ operator norm,
  because the dependence on $\lpr$ in \myeqref{eq:inoplpr}--\myeqref{eq:allop-def} is continuous.
  By \myeqref{eq:eps0} for $\eps>0$ sufficiently small there are no corner exponents in any corner with value $-\eps$ for any $\lpr\in\cli01$,
  so by \cite[theorem 4.1.2, p.\ 99]{nazarov-plamenevskii} all the operators on the path are Fredholm. Therefore their index is constant, in particular
  \begin{alignat}{9}  \find ( \allop(\textcolor{red}{1}) : \domsob\sobl{-\eps,-\eps} \rightarrow  \ransob\sobl{-\eps,-\eps} )  \topref{eq:find-dirneunegeps}{=}  +3  . \myeqlabel{eq:fredi-lpr}\end{alignat}
  Now we use \cite[theorem 4.3.3 (p.\ 112)]{nazarov-plamenevskii} to increase the corner weight $\wtoi$ in the $\pqo,\pqi$ corners from $-\eps$ to $1+\eps$:
  \begin{alignat}{9} \find ( \allop(1) : \domsob\sobl{\textcolor{red}{1+\eps},-\eps} \rightarrow  \ransob\sobl{\textcolor{red}{1+\eps},-\eps} )  \topref{eq:fredi-lpr}{=}  +3 -2\cdot 2 = -1 ,
  \myeqlabel{eq:find-onepluseps-minuseps} \end{alignat}
  because by \myeqref{eq:eps0} we cross a corner exponent $0$ of multiplicity $2$ in each of the two corners.
  
  2. 
  The domain is too small because exponents $1+\eps$ force $\pqpot$ and $\nabp\pqpot$ to be zero in $\pqo,\pqi$.
  Physically we would have $\xxn=\nabp\pqpot_{|\pq=\pqi}=0$,
  so the shock would not detach from the wedge.
  Clearly we need to permit some $\pqpot$ with nonzero value or gradient in $\pqo,\pqi$.
  
  To this end we recall the cutoffs $\cutcori$ above;
  in the $\cori=o,n$ corners we extend the domain by multiples of cut-off $\pqpotshift$ from \myeqref{eq:pqpotshift}:
  \begin{alignat}{9} \extdomsob\sobl{1+\eps,\wtw} \defeq \domsob\sobl{1+\eps,\wtw} + \RR\cutfni\pqpotshift + \RR\cutfno\pqpotshift . \myeqlabel{eq:extdom-def}\end{alignat}
  These sums are direct because $\cutcori\pqpotshift$ for $\cori=o,i$ has nonzero gradient in the $\cori$ corner, incompatible with the $\wtoi=1+\eps$ corner weight. 
  Unlike $\pqpotshift$ (Proposition \myref{prop:pqpotshift}) the cut-off $\cutcori\pqpotshift$ need not be global solutions,
  but they are local solutions near every corner $\pq^\cori$ where $\cutcori$ is constant $0$ or $1$ (our problem is linear homogeneous). 
  So they are mapped by the operator $\allop(1)$ into functions that are zero near each corner, in particular elements of the \emph{same} codomain $\ransob\sobl{1+\eps,-\eps}$.
  In contrast the domain did increase by $2$ dimensions, so the Fredholm index increases\footnote{%
  If we extend the domain $\mathcal X$ of an operator $\mathcal T$ by multiples of $u$, the kernel increases if and only if $\mathcal Tu\in \mathcal TX$, i.e.\ the range does not increase. Either way the index increases by $1$.}
  by $2$:
  \begin{alignat}{9}&
  \find \big( \allop(1) : \textcolor{red}{\mathcal{\hat X}}^\sobl_{1+\eps,-\eps} \rightarrow \ransob\sobl{1+\eps,-\eps} \big)
  &\topref{eq:find-onepluseps-minuseps}{=}
  -1
  +2
  =
  +1
  \end{alignat}
  Now that the index is positive, each kernel $N^\sobl$ of $\allop(1)$ on $\extdomsob\sobl{1+\eps,-\eps}$ has positive dimension.
  \cite{nazarov-plamenevskii} state theorems for finite Sobolev exponent $\sobl$,
  but we may argue that as $\sobl\conv\infty$ the domains and therefore kernels $N^\sobl$ are decreasing, so their intersection $N^\infty$ also has positive dimension;
  besides, kernels of Fredholm operators are finite, so for sufficiently large $\sobl$ we have $N^\sobl=N^\infty$.
  The dimension of the cokernels is irrelevant because our problem is homogeneous, with zero data. 
  Pick some nonzero
               
  \begin{alignat}{9} \pqpot' \in  N^\infty \subset \extdomsob\infty{1+\eps,-\eps} . \myeqlabel{eq:pqpot-prime-space} \end{alignat}
  This solution $\pqpot'$ is in particular $\Cinf$ on $\cDompq\setdiff\sett{\pqw,\pqi,\pqo}$,
  and by choice of corner weights $1+\eps$ in $\pqo,\pqi$ (and linear $\pqpotshift$) it is $\Conce(\cDompq\setdiff\sett{\pqw})$.
  As explained in Remark \myref{rem:analyticity} $\pqpot'$ is in fact analytic except in the corners since the boundary conditions and equation have analytic coefficient functions.

  3. 
  We apply \cite[Theorem 4.2.1 (p.\ 106)]{nazarov-plamenevskii} to obtain the leading terms in the weak corner $\pqw$.
  We first increase the corner weight $\wtw$ from $-\eps$ to $\eps$,
  crossing (by \myeqref{eq:eps0}, \myeqref{eq:pqpot-w-zeroplus} and \myeqref{eq:critex-w-zerominus})
  only the corner exponent $q^w_{0-}=0$ corresponding to the principal-part solution $\pqpot^w_{0-}=1$ \myeqref{eq:pqpot-w-zerominus}:
  for some scalar constant $c^w_{0-}$
  \begin{alignat}{9} 
  \cutfnw\cdot ( \subeq{ \pqpot' }{\mbox{$ \in\extdomsob{\infty}{1+\eps,\textcolor{red}{-\eps}} $}}
  -
  c^w_{0-} \subeq{\pqpot^w_{0-}}{=1} )
  &\in
  \extdomsob{\infty}{1+\eps,\textcolor{red}{\eps}} \quad.
  \myeqlabel{eq:stepover-weak-0-minusB}
  \end{alignat}
  Here and later their theorem is applied to the part of $\pqpot$ belonging to the $\domsob\sobl{1+\eps,-\eps}$ summand in \myeqref{eq:extdom-def},
  with other summands unchanged, using that they are local solutions near every corner;
  the theorem is stated for finite $\sobl$, but each yields the same $c$ because the $c$ terms have too little decay for the right-hand side spaces,
  so we may take the intersection over $\sobl$.
  
  To make progress we need to increase the corner weight further and find the next deeper asymptotic term, by applying the theorem again.
  But it applies only to local solutions of the \emph{full} problem,
  whereas the leading asymptotic term $c^w_{0-}\pqpot^w_{0-}$ solves the \emph{principal part}. 
  Although it is generally difficult to find explicit leading terms that are full solutions,
  here we have $\pqpotshift$ \myeqref{eq:pqpotshift}, with value $\pqwx\neq 0$ as $\pq\conv\pqw$ (note $\pqix=0$).
  The global full solution
  \begin{alignat}{9}  \pqpot''  &:=  \pqpot' - c^w_{0-} \frac{\pqpotshift}{\pqwx}      \myeqlabel{eq:stepover-firstline} \end{alignat}
  satisfies
  \begin{alignat}{9} \pqpot'' & \in  \extdomsob{\infty}{1+\eps,\textcolor{red}{\eps}} \end{alignat}
  because $\pqpotshift$ is linear, hence in every Sobolev space; $\pqpotshift$ near the corners $\pqo,\pqi$ is permitted by \myeqref{eq:extdom-def},
  as is $\pqpot'$ by \myeqref{eq:pqpot-prime-space} (same weight $1+\eps$ in those corners); 
  near $\pqw$ we argue
  \begin{alignat}{9} \cutfnw \pqpot''
  &\toprefb{eq:stepover-firstline}{eq:stepover-weak-0-minusB}{=}
  \subeq{  \cutfnw\cdot   (  \pqpot' - c^w_{0-} \pqpot^w_{0-}  )  }{\mbox{$
      \topref{eq:stepover-weak-0-minusB}{\in}  \extdomsob{\infty}{1+\eps,\textcolor{red}{\eps}}  $}}
  +  c^w_{0-}  \subeq{  \cutfnw\cdot   (  \supeq{ \supeq{\pqpot^w_{0-}}{=1} - \supeq{ \pqpotshift / \pqwx}{\rightarrow -1} ) }{=O(r)}   }{\mbox{$
      \in  \extdomsob{\infty}{1+\eps,\textcolor{red}{\eps}}  $}}.
  \myeqlabel{eq:stepover-weak-0-minus-reduced}
  \end{alignat}

  4. 
  Now we may apply \cite[Theorem 4.2.1]{nazarov-plamenevskii} repeatedly to increase the weak corner weight from $+\eps$ to $1+\eps$,
  crossing (by \myeqref{eq:eps0}, \myeqref{eq:pqpot-w-zeroplus} and \myeqref{eq:critex-w-zerominus}) the single exponent $q^w_{0+}\in\boi01$
  corresponding to the principal-part solution $\pqpot^w_{0+}$ \myeqref{eq:pqpot-w-zeroplus}: for some scalar constant $c^w_{0+}$
  \begin{alignat}{9} 
  \cutfnw\cdot (  \subeq{ \pqpot'' }{\mbox{$ \in \extdomsob{\infty}{1+\eps,\textcolor{red}{\eps}} $}}
  - c^w_{0+} \pqpot^w_{0+} )
  &\in 
  \extdomsob{\infty}{1+\eps,\textcolor{red}{1+\eps}}
  \myeqlabel{eq:xpotprimeprime-space}
  \end{alignat}
  If $c^w_{0+}=0$, then this 
  shows by definition \myeqref{eq:extdom-def} that
  $\pqpot'' \in \Conce(\cDompq) $, 
  in particular at the infinity corner $\pqw$, where $\pd2\pqpot''$ has a finite limit;
  if $c^w_{0+}\neq 0$, then \myeqref{eq:yv-particular-blowup} shows the limit exists but is infinite. 

  5. 
  By \myeqref{eq:xpotprimeprime-space} and definition \myeqref{eq:extdom-def} we have
  \begin{alignat}{9} \nabp\pqpot'' = c_o \subeq{\nabp\pqpotshift}{=(1,0)} + o(1) 
  \quad\text{as $\pq\conv\pqo$} \end{alignat}
  for some scalar constant $c_o$, but we want zero gradient in that corner, so we need to eliminate $c_o\pqpotshift$. For this we define
  \begin{alignat}{9} \pqpot'''  \defeq  \pqpot'' - c_o\pqpotshift , \myeqlabel{eq:pqpot-tripleprime}\end{alignat}
  which is still a global full solution because $\pqpot'',\pqpotshift$ are, but satisfies $\nabp\pqpot'''=0$ at $\pq\conv\pqo$,
  and because $\pqpotshift$ is linear we retain the well-defined limit of $\pd2\pqpot'''$ at $\pqw$, and that $\pqpot'''$ is $\Conce(\cDompq)$ if $\pqpot''$ is; besides
  \begin{alignat}{9} \pqpot''' \in \extdomsob\infty{1+\eps,-\eps}. \myeqlabel{eq:tripleprime-space} \end{alignat}

  6. 
  Finally we take $\pqpot=\pm\pqpot'''$, with $\pm$ chosen to make the limit $\yew$ of $\pd2\pqpot$ at $\pq\conv\pqw$ nonnegative (allowing $+\infty$), showing that either \myeqref{eq:eitheryblowup} or \myeqref{eq:orynonnegative} holds, and still $\pd1\pqpot=0$ at $\pq\conv\pqo$ showing \myeqref{eq:pqpot-nabp-zero-pqo}.

  By \myeqref{eq:tripleprime-space} $\pqpot\in\extdomsob\infty{1+\eps,-\eps}$, so
  \begin{alignat}{9} \pd2\pqpot \conv 0 \quad\text{as $\pq\conv\pqo,\pqi$;} \myeqlabel{eq:yv-pqi-pqo} \end{alignat}
  this holds for the $\domsob\infty{1+\eps,-\eps}$ part due to $\wtoi=1+\eps$, while the $\pqpotshift$ have $\pd2\pqpotshift=0$ everywhere \myeqref{eq:pqpotshift}.

  On $\cAxispq$ we have $\pqvar^1=0$, so the slip condition \myeqref{eq:bcond-Axis-scaled} reduces to
  \begin{alignat}{9} 0 = \pqvar^2\pd2\pqpot - \pqpot \end{alignat}
  which has solution $\pqpot=C\pqvar^2$ for some constant $C$. 
  By \myeqref{eq:yv-pqi-pqo} 
  only $C=0$ is possible. Thus \myeqref{eq:pqpot-zero-Axis} has also been shown and the proof is complete. 
\end{proof}

\mysection{Injectivity}\mylabel{section:pq-to-xx}

It is generally difficult to prove flows are univalent, or that solution constructed by existence proofs in the velocity or mass flux plane permit injective maps back to the $\xx$ plane.
But we can introduce a rather powerful way of proving injectivity using the device of collapsing level sets into edges, with applications far beyond this paper. 

From here till the end of the paper fix $\pqpot$ to be one of the solutions of Proposition \myref{prop:pqexi}.

\begin{proposition}
  \begin{alignat}{9} \yew>0  \myeqlabel{eq:yvw-nonzero} \end{alignat}
  and 
  \begin{alignat}{9} \ye = \pd2\pqpot \in \boi{0}{\yew} \quad\text{in $\cDompq\setdiff(\cAxispq\union\sett{\pqw})$.} \myeqlabel{eq:yv-positive-notonAxis}\end{alignat}
\end{proposition}\begin{proof}
  By Proposition \myref{prop:dex-vertical-ok}(iii) $\ye$ can only attain local extrema at $\cAxispq$, where it is zero, or in $\pqw$ where it has limit $\yew$.
  So if $\yew=0$, then $\ye=\pd2\pqpot=0$ in the entire domain, quickly leading to contradictions:
  \begin{alignat}{9}
  \pd 1\pqpot &\topref{eq:pqpot-nabp-zero-pqo}{=} 0 \topref{eq:pqpot-zero-Axis}{=} \pqpot \quad\text{in $\pq=\pqo$}\csep
  0 \topref{eq:interior-pq}{=} \adjt\Aco:\nabp^2\pqpot = \supeq{\Acosym^{22}}{>0}\pd 1^2\pqpot
  \\\impl\quad
  \pqpot &= 0.
  \end{alignat}
  So $\yew\neq 0$, and $\yew\geq 0$ by \myeqref{eq:eitheryblowup}+\myeqref{eq:orynonnegative} shows $\yew>0$,
  then Proposition \myref{prop:dex-vertical-ok}(iii) shows \myeqref{eq:yv-positive-notonAxis}.
\end{proof}

\begin{proposition}
  \mylabel{prop:singlecurve}%
  (i) For $\ye_+\in\roi{0}{\yew}$ each of the level sets
  \begin{alignat}{9} L := \selset{\pq\in\cDompq\setdiff\sett{\pqw}}{\pd 2\pqpot(\pq) = \ye_+}. \end{alignat} 
  is a simple curve that is analytic including at the endpoints, one of which is on $\cWedgepq\setdiff\sett{\pqw}$, the other on $\cShockpq\setdiff\sett{\pqw}$. \\
  (ii) $\ye=\pd 2\pqpot$ is strictly increasing along $\Shockpq$ from $\pqi$ to $\pqw$, and along $\Wedgepq$ from $\pqo$ to $\pqw$,
  with values from $0$ to $\yew$. 
\end{proposition}
\begin{proof}
  Consider $\ye_+\in\boi{0}{\yew}$. 
  Having derived the elliptic PDE \myeqref{eq:dex-nabp-pqpot-pde} for $\ye$ 
  we can apply Proposition \myref{prop:levels-tree} again:
  if one of the trees of $L$ has a vertex, necessarily of degree at least four, then the tree has at least four simple paths to the boundary with pairwise disjoint interiors. 
  On $\cAxispq$ we have $\ye=0<\ye_+$ \myeqref{eq:pqpot-zero-Axis}, and at $\pqw$ we have $\ye\conv\ye_w>\ye_+$, 
  so since $\Wedgepq,\Shockpq$ are analytic every path approaches a single endpoint on one of them (Remark \myref{rem:single-endpoint-boundary}). 
  Either $\Wedgepq$ or $\Shockpq$ is approached by at least two paths.

  If any tree has two or more paths to $\Wedgepq$, then we can form a simple loop $\bdry K$ with the segment of $\Wedgepq$ between the two endpoints (not necessarily distinct),
  enclosing by the Jordan curve theorem some nonempty set with compact closure $K$, and $\ye$ attains a maximum and a minimum over $K$.
  If both are on $L$ then $\ye$ is constant $\ye_+$ on $K$, so either way there is a local extremum in $K\setdiff L$, 
  but by Proposition \myref{prop:dex-vertical-ok}(iii) that is not possible.

  The same argument rules out two paths to $\Shockpq$.
  So there are no vertices; each $L$ tree is a single edge.
  By the same argument we cannot have edges running from $\Shockpq$ back to $\Shockpq$, or from $\Wedgepq$ back to $\Wedgepq$,
  so every edge runs from $\Shockpq$ to $\Wedgepq$. 
  If there are multiple edges with disjoint interior, then adding the $\Wedgepq$ segment between their endpoints there, as well as the $\Shockpq$ segment between those endpoints,
  we obtain another simple loop $\bdry K$ and another contradiction.

  For $\ye_+=0$ the level curve is exactly $\Axispq$, by \myeqref{eq:pqpot-zero-Axis} and \myeqref{eq:yv-positive-notonAxis}. 

  We have shown every $\ye\in\boi{0}{\yew}$ occurs exactly once on $\Shockpq$ and exactly once on $\Wedgepq$, proving $\ye$ is strictly increasing along each.
\end{proof}

\begin{proposition}
  \mylabel{prop:det-nonzero-xmonoylevel}%
  (i) $\Det\nabp^2\pqpot < 0$ on $\cDompq\setdiff\sett{\pqo,\pqi,\pqw}$.\\
  (ii) $\xe=\pd 1\pqpot$ is strictly increasing from $\Shockpq$ to $\Wedgepq$ along each $\ye$-level curve from Proposition \myref{prop:singlecurve} .\\
  (iii) $\pq\mapsto(\xe,\ye)$ is injective on $\cDompq\setdiff\sett{\pqw}$. 
\end{proposition}
\begin{proof}
  Obviously $\ye=\pd2\pqpot$ has along each $\ye$-level curve zero tangential derivative. 
  If in some point $\pq_*$ the tangential derivative of $\xe=\pd 1\pqpot$ is also zero, 
  or if $\Det\nabp^2\pqpot=0$ there, 
  then Lemma \myref{lemma:two-zero-three} and the interior pde \myeqref{eq:interior-pq} show $\nabp^2\pqpot=0$ in $\pq_*$,
  so $\nabp\ye=0$, with $\nabp^m\ye\neq 0$ for some minimal $m\geq 2$. 
  
  If $\pq_*$ is an interior point, then Proposition \myref{prop:levels-atpoint}(i) shows it is a vertex, 
  contradicting Proposition \myref{prop:singlecurve}.

  If $\pq_*$ is on $\Axispq$, where $\ye=0$, 
  then Proposition \myref{prop:levels-atpoint}(ii) shows at least $m-1\geq 1$ edges of level $0$ go to the interior $\Dompq$,
  contradicting \myeqref{eq:yv-positive-notonAxis}. 
  
  If $\pq_*$ is on $\Shockpq$ or $\Wedgepq$,
  then the boundary conditions \myeqref{eq:bcond-Shock} and \myeqref{eq:bcond-Wedge} have derivative direction $\pq-\pqi$ resp.\ $\pq$ nowhere tangential,
  by Proposition \myref{prop:dex-vertical-ok}(i)+(ii),
  so Proposition \myref{prop:levels-atpoint}(iii) shows that $m\geq 2$ edges go to the interior,
  again contradicting Proposition \myref{prop:singlecurve}.

  We have shown $\Det\nabp^2\pqpot\neq 0$ everywhere.
  But this is necessarily $<$, because $>0$ means $\nabp^2\pqpot$ is a definite symmetric matrix,
  but so is $\adjt\Aco$, violating the interior pde \myeqref{eq:interior-pq} $0 = \adjt\Aco:\nabp^2\pqpot$ (which holds by limit in smooth boundary points as well).
  We have proven part (i).

  Besides we have shown $\xe$ is strictly monotone along every $\ye$-level set for value $\ye_+\in\roi{0}{\yew}$, but those fill $\cDompq\setdiff\sett{\pqw}$,
  so $\pq\mapsto(\xe,\ye)$ is injective there, showing (iii).

  Since $\ye=\pd2\pqpot=0$ on $\Axispq$, along which a $\pd2$ tangential derivative shows $\pd2^2\pqpot=0$, necessarily $0\neq\pd1\pd2\pqpot$ (else $\Det\nabp^2\pqpot=0$).
  But then $\ye>0$ in $\cDompq\setdiff\cAxispq$ only permits $0>\pd1\ye=\pd1\pd2\pqpot=\pd2\xe$, so $\xe$ is strictly \emph{decreasing} from $\pqo$ to $\pqi$ along $\Axispq$.
  By continuity it must also be strictly decreasing along every other level set from $\Wedgepq$ to $\Shockpq$,
  because we already knew it is strictly monotone. 
  This completes part (ii).
\end{proof}

\begin{proof}[Proof of Theorem \myref{th:weak-inf-existence}]
  By Proposition \myref{prop:det-nonzero-xmonoylevel}(iii) the map $\pq\mapsto\xx=\nabp\pqpot$ is injective on $\cDompq\setdiff\sett{\pqw}$.
  Let $\Domxy$ be the image of $\Dompq$ under this map,
  with $\Shockxy,\Axisxy,\Wedgexy$ the images of $\Shockpq,\Axispq,\Wedgepq$.
  By Proposition \myref{prop:det-nonzero-xmonoylevel}(ii) $\Axisxy$ is a segment of the negative $\xe$ axis,
  because by \myeqref{eq:pqpot-nabp-zero-pqo} and \myeqref{eq:pqpot-zero-Axis} the image of the first endpoint $\pq=\pqo$ is $\xx=\xxo$;
  let $\xxn$ be the image of the other endpoint $\pqi$. 

  The map $\pq\mapsto\xx$ is analytic away from the corners like $\pqpot$ (Proposition \myref{prop:pqexi}),
  so by the analytic inverse function theorem \cite[theorem\ 2.5.1]{krantz-parks} the inverse is also analytic. 
  The map is $\Conce$ at the corners $\pqo,\pqi$, so the inverse is also $\Conce$ at $\xxo,\xxn$.
  By Proposition \myref{prop:det-nonzero-xmonoylevel}(ii) $\ye$ in strictly increasing along $\Shockpq,\Wedgepq$,
  so they are graphs of some $\Conce$ functions $\hs,\hw$ that are analytic away from the endpoints. 
  
  As in \myeqref{eq:stf-pqpot} define 
  \begin{alignat}{9} \stf(\xx) = \pq(\xx)\dotp\xx - \pqpot(\pq(\xx))   .  \end{alignat}
  As discussed in section \myref{section:pq-equations}, $\stf$ solves the original $\xx$-plane equation \myeqref{eq:interior-stf}
  in $\Domxy$ since $\pqpot$ solves the $\pq$-plane version \myeqref{eq:interior-pq},
  and it satisfies the slip condition \myeqref{eq:slip-stf} on $\Axisxy$ and $\Wedgexy$ because $\pqpot$ satisfies the analogue \myeqref{eq:slipcond-pq};
  the choice of $\Shockpq$ as a segment of the $\pq$ shock polar, combined with the boundary condition \myeqref{eq:polarcond} for $\pqpot$ implies
  that $\stf$ satisfies the shock conditions \myeqref{eq:shock-jmp-stf} (which implies \myeqref{eq:shock-jn}) and \myeqref{eq:shock-ut}.

  Assume that $\yew<+\infty$.
  Then by \myeqref{eq:orynonnegative} we have $\pqpot\in\Conce(\cDompq)$, including the $\pqw$ corner,
  so $\stf\in\Conce(\cDomxy)$ including the $\xxw$ corner that is the image of $\pqw$ under $\pq\mapsto\xx=\nabp\pqpot$;
  $\xxw$ is a common endpoint of $\Shockxy$ and $\Wedgexy$. 
  But we have already proven the shock relations which imply $\jj\dotp\shockn>0$ everywhere on $\Shockxy$, 
  while the slip condition implies $\jj\dotp\walln=0$ on $\Wedgexy$ and $\Axisxy$.
  Therefore $\int_{\bdry\Domxy}\jj\dotp\vec nds=\int_{\Domxy}\nabx\dotp\jj dA$ is nonzero so that conservation of mass $0=\nabla\dotp\jj$ is violated.
  This contradiction proves $\yew=+\infty$.
  Therefore $\Shockxy,\Wedgexy$ converge to infinity without meeting and $\Domxy$ is unbounded, with $\vv\conv\vvw$ at upper infinity by choice of corner $\pqw$. 
\end{proof}

\mysection{Appendix: Shock analyticity}\mylabel{section:shock-analyticity}

\begin{proposition}
  \mylabel{prop:shock-interior-analytic}%
  $\eShockxy$ is an analytic curve, and $\stf$ is analytic in $\eDomxy\union\eShockxy$. 
\end{proposition}  \begin{proof}  
  \newcommand{\avec}{{\vec a}}
  \newcommand{\naba}{\nabla_\avec}
  The shock is a free boundary, but it is standard to reduce such problems to fixed boundaries by some coordinate change.
  A basic idea is to find a quantity that has a fixed value everywhere on the shock by using some combination of the shock relations, for example \myeqref{eq:shock-jmp-vpot} $\vpoti-\vpot=\const$, or \myeqref{eq:shock-jmp-stf} $\stfi-\stf=0$, or any kind of combination thereof; see \cite[section 3.1]{elling-sonic-potf} or \cite{elling-liu-pmeyer} for other examples.
  Among the infinity of possible choices the following coordinates $\avec=(a,b)$ seem relatively convenient:
  \begin{alignat}{9} a &\defeq \vpoti-\vpot ,
  \\ b &\defeq \ye . \end{alignat}
  The Jacobian is
  \begin{alignat}{9} \avec_\xx
  =
  \begin{bmatrix}
    a_\xe & a_\ye \\
    b_\xe & b_\ye
  \end{bmatrix}
  &=
  \begin{bmatrix}
    \vpot_{i\xe}-\vpot_\xe & \vpot_{i\ye}-\vpot_\ye \\
    0 & 1
  \end{bmatrix}
  =
  \begin{bmatrix}
    \vxi-\vx & -\vy \\
    0 & 1
  \end{bmatrix}
  \myeqlabel{eq:ab-xy}\end{alignat}
  which is obviously invertible since $\vxi=|\vvi|>|\vv|\geq\vx$ for any admissible shock. Besides 
  \begin{alignat}{9}
  b_\ye &= 1 = \ye_b
  \csep
  \ye_a = 0 = b_\xe
  \csep
  \begin{bmatrix}
    a_\xe & a_\ye \\
    b_\xe & b_\ye
  \end{bmatrix}
  &=
  \begin{bmatrix}
    \xe_a^{-1} & -\xe_a^{-1}\xe_b \\
    0 & 1
  \end{bmatrix} ,
  \myeqlabel{eq:axmat}
  \end{alignat}
  so we can solve \myeqref{eq:ab-xy} for
  \begin{alignat}{9} \vx &= \vxi - a_\xe = \vxi - \xe_a^{-1} , \notag
  \\ \vy &= -a_\ye = \xe_a^{-1}\xe_b . \myeqlabel{eq:nabla-phi-xab}\end{alignat}
  So $\dens$, which is an analytic function of $\vv$ (section \myref{section:equations}), becomes an analytic function of $\naba\xe$. The continuity equation transforms to
  \begin{alignat}{9} 0
  &= \nabx\dotp( \dens \vv)  \myeqlabel{eq:ndiv-dens-vv}
  \\&= \naba\dotp(\dens\adjt\xx_\avec \vv) 
  \\&= \begin{bmatrix}
    \partial_a \\
    \partial_b
  \end{bmatrix}\dotp(\dens
  \begin{bmatrix}
    \ye_b & -\xe_b \\
    -\ye_a & \xe_a
  \end{bmatrix}
  \begin{bmatrix}
    \vx \\
    \vy
  \end{bmatrix}) 
  \\&\toprefb{eq:axmat}{eq:nabla-phi-xab}{=} 
  \begin{bmatrix}
    \partial_a \\
    \partial_b
  \end{bmatrix}
  \dotp(\dens\begin{bmatrix}
  1 & -\xe_b \\
  0 & \xe_a
  \end{bmatrix} \begin{bmatrix}
    \vxi - \xe_a^{-1}  \\
    \xe_a^{-1}\xe_b
  \end{bmatrix}) 
  \\&= \Big( \dens \big( \vxi - \frac{1 + \xe_b^2}{\xe_a}\big)\Big)_a + (\dens \xe_b)_b.
  \end{alignat}
  This is an analytic 2nd order equation for $\xe=\xe(a,b)$, quasilinear divergence-form. Since the non-divergence form $\Aco:\nabx^2\vpot=0$ of \myeqref{eq:ndiv-dens-vv}
  has symmetric positive definite $\Aco$,
  and $\xx\mapsto\avec$ is in particular a local $\Conce$ diffeomorphism,
  necessarily the new non-divergence form is still uniformly elliptic
  because its 2nd order part $(\avec_\xx\Aco\avec_\xx^T):\naba^2\vpot$ has symmetric positive definite coefficient matrix $\avec_\xx\Aco\avec_\xx^T$. 

  $\vpot-\vpoti=\const=a$, which is equivalent to $0\topref{eq:shock-ut}{=}\shockn^\perp\dotp\jmp{\nabx\vpot}$, was used to define $a$,
  so we use the other shock relation \myeqref{eq:shock-jn} as new boundary condition.
  $a$ is constant on the shock, so $\partial_b$ is a tangential derivative, and $(\xe_b,\ye_b)=(\xe_b,1)$ is an $\xx$-space shock tangent, so $(1,-\xe_b)$ is a normal:
  \begin{alignat}{9}  0  \topref{eq:shock-jn}{=}  (1,-\xe_b)\dotp\jmp{\dens\vv}  \myeqlabel{eq:x-bcond}  \end{alignat}
  with $\vv$ as in \myeqref{eq:nabla-phi-xab}. 

  As we show now this boundary condition is nowhere degenerate oblique.
  (This property not indispensable in two dimensions,
  but much of the literature assumes the condition to treat all dimensions at once without distinction, so it is convenient to have the property available to cite wider range of texts.
  Besides, the reduction above also works in higher dimensions, if $b$ is extended to $b_j=y_j$ for $j=2,...,n$, with obvious modifications, so the property \emph{is} needed there.)
  To check degenerate obliqueness linearize the boundary condition, with $\delta$ indicating the first variation of $x$ in $(a,b)$ coordinates:
  \begin{alignat}{9} \delta(\dens\vv-\densi\vvi)
  &=
  \delta ( \dens \vv ) - \supeq{ \delta(\densi\vvi) }{=0}
  \\&\topref{eq:ddensvv-dvv}{=}
  \dens \Aco \delta\vv
  \\&\topref{eq:nabla-phi-xab}{=}
  \dens \Aco \delta(\vxi-\xe_a^{-1},\xe_a^{-1}\xe_b)
  \\&=
  \subeq{\xe_a^{-2}}{>0} \subeq{\dens}{>0} \Aco (1,-\xe_b) \delta\xe_a + ... 
  \end{alignat}
  where $...$ are terms without $\delta\xe_a$ (i.e.\ coefficients times $\delta\xe_b$). $\delta(1,-\xe_b)=...$ has no $\delta\xe_a$, so the entire right-hand side of \myeqref{eq:x-bcond} linearizes to
  \begin{alignat}{9} \xe_a^{-2}\dens~\subeq{(1,-\xe_b)\dotp\Aco(1,-\xe_b)}{>0} \delta\xe_a + ... ~.\end{alignat}
  Hence the coefficient of $\delta\xe_a$ is always nonzero, so the nonlinear boundary condition \myeqref{eq:x-bcond} is not degenerate oblique, in fact uniformly not, by smoothness of the coefficient functions. 

  Having reduced the problem to a fixed boundary, we may apply standard results, such as \cite[section 6.7]{morrey-book} (which would permit degenerate oblique conditions, because in two dimensions they satisfy the ``complementing condition'', also known as ``Lopatinski-Shapiro condition'' among many other names). 
  This shows $\xe$ is analytic in a neighbourhood of any $(0,b)$, so $\avec\mapsto\xx$ is analytic with analytic inverse, the image $\Shockxy$ of $\sett{a=0}$ is locally a nonsingular analytic curve, $\vpot_\xx(\xx)=\vpot_\avec(\avec(\xx))\avec_\xx(\xx)$ is also analytic, and so is $\vpot$.
  (We could also invoke literature such as \cite[Theorem 6.26]{gilbarg-trudinger}, which does require uniform non-degeneracy,
  to improve regularity from various $C^{m,\alpha}$ if we had not assumed $\Cinf$.)
\end{proof}

\section*{Data availability statement}
This is a purely mathematical work and no data was created or analyzed in this study.
\section*{Acknowledgements}
This research was partially supported by Taiwan NSTC grant 113-2115-M-001-005-MY3.

\bibliographystyle{unsrt}

%\bibliography{../../../pmeyer/elling.bib}

\end{document}